\documentclass[twoside,11pt]{article}

\usepackage[usenames]{color} 
\definecolor{Plum}{rgb}{.4,0,.4} 
\definecolor{BrickRed}{rgb}{0.6,0,0} 
\RequirePackage{amsmath,amsfonts,amssymb}
\RequirePackage{graphicx}
\RequirePackage{natbib} 
\numberwithin{equation}{section}

\usepackage{bbm}
\usepackage{tabularx,multirow,booktabs}
\usepackage{caption}

\DeclareMathOperator*{\argmin}{arg\,min}

\DeclareMathOperator*{\E}{\mathbb{E}}
\newcommand{\F}{\mathcal{F}}
\newcommand{\D}{\mathcal{D}}
\newcommand{\G}{\mathcal{G}}
\newcommand{\cH}{\mathcal{H}}
\newcommand{\cR}{\mathcal{R}}
\newcommand{\cT}{\mathcal{T}}
\newcommand{\bW}{\mathbf{W}}

\DeclareMathOperator{\sign}{sign}

\usepackage{jmlr2e}

\usepackage{lastpage}
\jmlrheading{22}{2021}{1-\pageref{LastPage}}{8/20; Revised
4/21}{9/21}{20-911}{Tengyuan Liang}

\ShortHeadings{Statistical Theory for GANs}{Liang}
\firstpageno{1}

\begin{document}

\title{How Well Generative Adversarial Networks Learn Distributions}

\author{\name Tengyuan Liang \email tengyuan.liang@chicagobooth.edu \\
       \addr Econometrics and Statistics\\
       University of Chicago, Booth School of Business\\
       Chicago, IL 60637, USA}

\editor{Ambuj Tewari}

\maketitle

\begin{abstract}
This paper studies the rates of convergence for learning distributions implicitly with the adversarial framework and Generative Adversarial Networks (GANs), which subsume Wasserstein, Sobolev, MMD GAN, and Generalized/Simulated Method of Moments (GMM/SMM) as special cases. We study a wide range of parametric and nonparametric target distributions under a host of objective evaluation metrics. We investigate how to obtain valid statistical guarantees for GANs through the lens of regularization. On the nonparametric end, we derive the optimal minimax rates for distribution estimation under the adversarial framework. On the parametric end, we establish a theory for general neural network classes (including deep leaky ReLU networks) that characterizes the interplay on the choice of generator and discriminator pair. We discover and isolate a new notion of regularization, called the generator-discriminator-pair regularization, that sheds light on the advantage of GANs compared to classical parametric and nonparametric approaches for explicit distribution estimation. We develop novel oracle inequalities as the main technical tools for analyzing GANs, which are of independent interest.
\end{abstract}


\section{Introduction}

Generative models such as Generative Adversarial Networks \citep{goodfellow2014generative, li2015generative, arjovsky2017wasserstein, dziugaite2015training} have recently stood out as an important unsupervised method for learning and efficient sampling from a complex target data distribution. Despite the celebrated empirical success, many questions on the theory \citep{liu2017approximation, liang2017well, singh2018nonparametric, liu2018inductive} and mechanism of GANs \citep{arora2017gans,arora2017generalization,daskalakis2017training, mescheder2017numerics} remain to be elucidated.

At the population level, one general formulation of the adversarial framework \citep{arjovsky2017wasserstein,li2015generative, dziugaite2015training, liu2017approximation, mroueh2017sobolev} considers the following minimax problem,
\begin{align*}
	\min_{\mu \in \D_{G}} \max_{f \in \F_{D}} \E_{Y \sim \mu} f(Y) - \E_{X \sim \nu} f(X).
\end{align*}
In plain language, given a target probability distribution $\nu$, one seeks a simulated probability distribution $\mu$ from a \textit{generator class} $\D_G$, so to minimize the loss incurred by a host of test functions inside a \textit{discriminator class} $\F_D$. In practice, both \textit{the generator and the discriminator classes} are parametrized by deep neural networks. To be concrete, $\D_G$ quantifies the implicit distributions realized by neural network transformations that push forward simple input random variables, for instance, with multi-dimensional uniform or Gaussian distribution. $\F_D$ represents certain neural network functions. In practice, one only has access to finite samples of the target distribution $\nu$. Let us denote $\widehat{\nu}^n$ as the empirical distribution based on $n$ i.i.d. samples from $\nu$, then the adversarial framework solves the following empirical problem
\begin{align}
	\label{eq:gan}
	\widehat{\mu} \in \argmin_{\mu \in \D_{G}} \max_{f \in \F_{D}} \E_{Y \sim \mu} f(Y) - \E_{X \sim \widehat{\nu}^n} f(X).
\end{align}
Here the first expectation over $Y\sim \mu$ can be calculated efficiently using simulations with arbitrary accuracy, since samples from $\mu$ can be simulated directly by pushing-forward random inputs. 
A natural question is to understand, how well the simulated distribution $\widehat \mu$ estimates the target $\nu$, under a host of evaluation metrics.

In machine learning language, the adversarial framework is induced by a certain Integral Probability Metric (IPM) quantifying the closeness between probability distributions. Define the IPM for a symmetric function class $\F$ as
\begin{align*}
	d_{\F}(\mu, \nu) := \sup_{f \in \F} \E_{Y \sim \mu} f(Y) - \E_{X \sim \nu} f(X)  =  \sup_{f \in \F} \int_{\Omega}  f (d\mu - d\nu).
\end{align*}
By choosing different $\F$'s, the adversarial framework can express a host of commonly-used metrics. To name a few, (1) Wasserstein GAN \citep{arjovsky2017wasserstein}: $\mathcal{F}$ consists of Lipschitz-$1$ functions, and the IPM is the Wasserstein-1 metric $d_W(\cdot, \cdot)$. (2) Maximum Mean Discrepancy (MMD) GAN \citep{dziugaite2015training,li2015generative, arbel2018gradient}: let $\mathcal{H}$ be a Reproducing Kernel Hilbert Space (RKHS), and $\mathcal{F}$ consists of functions with bounded RKHS norm $\F = \{ f \in \cH ~|~ \| f \|_{\cH} \leq 1 \}$. (3) Sobolev GAN \citep{mroueh2017sobolev}: $\mathcal{F}$ is the Sobolev class with certain smoothness. (4) Total Variation metric $d_{TV}(\cdot, \cdot)$: $\mathcal{F}$ represents all functions bounded by $1$. 
We refer the readers to \cite{liu2017approximation} for other related formulations of GAN.
Conceptually, the discriminator function class induces a collection of ``moment conditions'' assessing the closeness between distributions, as in the Generalized Method of Moments (GMM) \citep{hansen1982LargeSample}.

In the statistical literature, explicit distribution estimation, or density estimation, has been a fundamental topic in nonparametric statistics \citep{nemirovski2000topics, tsybakov2009introduction, wassermann2006all} and in parametric models \citep{brown1986fundamentals}. In the parametric case, learning density simply reduces to parameter estimation. In the nonparametric case, the optimal minimax rates have been established for a wide range of density function classes quantified by the smoothness property \citep{stone1982optimal}. However, it is not practical to simulate samples efficiently from these minimax optimal explicit density estimators, especially for multi-dimensional data. 

The econometrics literature has explored an alternative implicit distribution estimation approach \citep{back1993ImpliedProbabilities,imbens1995information} using the Simulated Method of Moments (SMM) \citep{pakes1989SimulationAsymptotics, mcfadden1989MethodSimulated}. Such an SMM approach turns out to be a special case of GANs in formulation \eqref{eq:gan}: SMM implicitly estimates the target distribution with a simulated distribution (from a certain parametric class $\mu_{\widehat{\theta}} \in \D_{G}$), by matching the moment conditions (induced by functions $f \in \F_{D}$) to the empirical distribution $\widehat{\nu}^n$. In the classic method of moments with finite $K$ moment conditions $\{ \phi_k(x), k\in [K] \}$, $\F_{D}$ consists of functions satisfying a quadratic constraint $\{ f(x) = \sum_{k\in [K]} \omega_k \phi_k(x) ~|~ \omega^\top \bW^{-1} \omega \leq 1, \omega \in \mathbb{R}^K \}$ with a given symmetric positive-definite weight matrix $\bW \in \mathbb{R}^{K \times K}$. In this language, the adversarial framework in \eqref{eq:gan} extends the SMM to where the moment conditions are induced by a rich class of functions $\F_{D}$. More recently,  \cite{athey2019using} conducted a systematic empirical study to learn the distributions of real economic data sets using Wasserstein GAN, suggesting the effectiveness of such an implicit distribution estimation approach in modern practice, even beyond the typical computer vision domain.

The current paper studies the \textit{Adversarial Framework} and \textit{Generative Adversarial Networks} for implicitly learning distributions from a statistical vantage point. As discussed in the paragraphs before, the problem is fundamental to statistics, machine learning, and econometrics. 
We intend to answer the following questions: 
\begin{enumerate}
	\item How well do GANs learn a wide range of target distributions (both in nonparametric and parametric cases), under a collection of objective evaluation metrics?
	\item How to leverage the adversarial framework to achieve better theoretical guarantees through the lens of regularization? 
\end{enumerate}
We discover and isolate a new notion of regularization, called the \textit{generator-discriminator-pair regularization}, which provides rigorous guidance on balancing the complexities of the generator and discriminator. We emphasize that several curious features of this pair regularization appear to be new to the literature. As a unified theme in theory, we develop oracle inequalities for analyzing the generative adversarial framework, which could be of independent interest for further theoretical research on GANs. It is worth noting that the early draft of this paper formulated the first statistical framework to study GANs.

\subsection{Contributions and Organization}
\label{sec:contribution}

The paper is organized into two main parts: the \textit{Adversarial Framework} and the \textit{Generative Adversarial Networks}.

\paragraph{\it Roadmap of Results and Overall Goal}
Our goal is to provide a comprehensive statistical treatment of the adversarial framework and GANs under two important settings. First, the generator and discriminator fall under the nonparametric classes, studied under the adversarial framework. Second, the generator and discriminator are the classes parametrized by neural networks as in GANs. We summarize in Table~\ref{table:summary} a roadmap of results for readers to navigate.  We emphasize that all the theoretical results (including Wasserstein, Sobolev and MMD GAN, GMM and SMM) follow a unified oracle inequality approach, demonstrating the universality of our framework.
In Table~\ref{table:summary}, we reserve the following symbols for the following properties of the theorems.
\begin{align}
	\label{eq:special-symbols}
	& ({\color{red} \G \dagger}): \quad \text{generator $\G$ could be mis-specified for $\nu$, $\nu \notin \G$.} \\
	& ({\color{blue} \F \ddagger}): \quad \text{discriminator $\F$ could be mis-specified for the metric, $d_{\F} \neq d_{eval}$.} \nonumber\\
	& ({\color{Plum} m \ast}): \quad \text{the result accounts for finite $m$ samples of the generator.}  \nonumber
\end{align} 

The main technical contributions are the development of the oracle inequalities for analyzing GANs, and the formulation of the new generator-discriminator-pair regularization. 

\begin{table}[ht!]
	\captionsetup{font=scriptsize}
	\scriptsize
	\newcolumntype{Y}{>{\centering\arraybackslash}X}
    \begin{tabularx}{\columnwidth}{@{} m{2cm} Y  Y Y| Y Y c @{}}
        \toprule
        Goal & Evaluation Metric & \multicolumn{2}{c|}{Results}  & Generator Class $\G$ & Discriminator Class $\F$ & ~~~Property  \\
        \midrule
		\multirow[m]{3}{=}{Adversarial Framework (nonparametric)} & \multirow{3}{*}{$d_{\F}$} 
		& Sobolev GAN & minimax optimal (Thm. \ref{thm:optimal-rates-sobolev})  & Sobolev $W^\alpha$ & Sobolev $W^{\beta}$ &  \\ 
	    \cmidrule{3-7}
		&  & MMD GAN & upper bound (Thm. \ref{thm:optimal-rates-rkhs}) & smooth subclass in RKHS & RKHS $\cH$ & \\
		\cmidrule{3-7}
		&  &  & oracle results (Thm. \ref{thm:nonparam-gans}) & any & Sobolev $W^{\beta}$ & $\color{red} \G \dagger$ \\
		\cmidrule{3-7}
		&  & GMM or SMM & oracle results (Cor. \ref{coro:GMM}) & any & moment conditions & $\color{red} \G \dagger$ \\
		\midrule
		\multirow[m]{3}{=}{Generative Adversarial Networks (parametric)}
		& $d_{TV}$ & leaky-ReLU GAN & upper bound (Thm.~\ref{thm:leaky-ReLU}) &  leaky-ReLU & leaky-ReLU &  $\color{blue} \F \ddagger$, $\color{Plum} m \ast$\\
		\cmidrule{2-7}
		& $d_{TV}, d_{JS}, d_{H}$  & any GAN or SMM & oracle results  (Thm.~\ref{thm:param-gans} \& \ref{thm:param-gans-hellinger}) &  general & general & $\color{red} \G \dagger$, $\color{blue} \F \ddagger$, $\color{Plum} m \ast$\\
		\cmidrule{2-7}
		& $d_{W}$ & Lipschitz GAN & oracle results (Cor.~\ref{cor:wass}) & Lipschitz neural networks & Lipschitz neural networks & $\color{red} \G \dagger$, $\color{blue} \F \ddagger$, $\color{Plum} m \ast$\\
        \bottomrule
    \end{tabularx}
 	\caption{Roadmap of results. The symbols are defined in \eqref{eq:special-symbols}: ($\color{red} \G \dagger$) and ($\color{blue} \F \ddagger$) denoting the mis-specification for the generator class and the discriminator class respectively, and ($\color{Plum} m \ast$) indicates the dependence on the number of simulated samples.}
	 \label{table:summary}
\end{table}

\paragraph{\it Adversarial Framework}
One key component of GAN is the adversarial framework: evaluating the performance of the learned distribution by the adversarial loss. 
Under the adversarial loss $d_{\F_D}(\cdot, \cdot)$ (IPM induced by the discriminator class $\F_D$), we study the minimax optimal rates for learning the target distribution $\nu$ based on $n$-i.i.d. samples. We formulate such an adversarial framework following the classic nonparametric literature by considering a host of nonparametric target distributions $\mu$ and discriminator classes $\F_D$ quantified by their smoothness property. Using an oracle inequality, we extend to the case when the generator class $\D_{G}$ misspecifies the target distribution $\nu$, for the procedure
\begin{align}
	\label{eq:adv-frame}
	\widehat{\mu}_n = \argmin_{\mu \in \D_{G}} \max_{f \in \F_{D}} \E_{Y \sim \mu} f(Y) - \E_{X \sim \widehat{\nu}^n} f(X).
\end{align} 

Our contributions are: (1) we derive the optimal minimax rates of the adversarial framework for learning a host of nonparametric distribution families, and how to attain the optimal rates; (2) we show how the regularities of target $\nu$ and of the class $\F_{D}$ affects the minimax rates explicitly, and under what cases fast rates are possible; (3) As a byproduct, we show that GMM and SMM can be derived as a particular case of the adversarial framework, and obtain explicit non-asymptotic rates for GMM.

\paragraph{\it Generative Adversarial Networks}
In practice, GANs are parametrized by deep neural networks.
Building upon the adversarial framework, we directly analyze the rates for the following parametrized GANs estimator with the generator network $\G$ (parametrized by $\theta$) and discriminator network $\F$ (parametrized by $\omega$)
	\begin{align}
		\label{eq:gan-network}
		\widehat{\theta}_{m,n} \in \argmin_{\theta: g_\theta \in \G} \max_{\omega: f_\omega \in \F} ~~ \left\{ \widehat{\E}_m  f_\omega (g_\theta(Z))  - \widehat{\E}_n  f_\omega(X)  \right\}.
	\end{align}
Here $m$ and $n$ denote the number of the simulated generator samples and target distribution samples. We emphasize two key facts about this procedure. First, the distribution estimator is implicit: the estimator is the probability distribution of a random variable push-forwarded by the transformation map $Z \mapsto g_{\widehat{\theta}_{m,n}}(Z)$ with $g_{\widehat{\theta}_{m,n}}: \mathcal{Z} \rightarrow \mathcal{X}$ and $Z \sim \pi$ (input uniform distribution). The theory for the implicit distribution estimator (such as GANs) is missing the current literature. Second, the objective evaluation metrics we investigate include Jensen-Shannon divergence $d_{JS}$, Total Variation $d_{TV}$, Wasserstein $d_W$, and Hellinger distances $d_{H}$, which all differ from and are misspecified by the generator $\F$. 

Our contributions are: (1) we derive the parametric rates on the closeness between the implicit distribution estimator (distribution of $g_{\widehat{\theta}_{m,n}}(Z)$) and the target $\nu$ under objective metrics, when both $\G$ and $\F$ are parametrized by general neural networks; (2) We rigorously formulate the complex trade-offs on the choices of the generator $\G$ and the discriminator $\F$ as a \textit{pair regularization}. We evaluate how this new notion of regularization affects the rates for GANs; 
(3) As a direct application of the general theoretical framework, we showcase how to identify good $(\G, \F)$ pairs to obtain fast parametric rates using two extreme examples: (a) learning distributions realizable by deep leaky ReLU networks, and (b) learning multivariate Gaussian distributions with two-layer networks. In both cases, the upper rates we obtain provide optimal sample complexity (up to logarithmic factors).

\paragraph{\it Organization} Finally, the paper is organized as follows. Section~\ref{sec:nonparam} consists of main nonparametric results and the adversarial framework. Section~\ref{sec:param} contains the main parametric results for GAN with neural network generator and discriminator classes, where we introduce the new notion of pair regularization. Further discussions on the generator-discriminator-pair regularization and connections to the regularity theory in optimal transport is deferred to Section~\ref{sec:discussion}. The main proofs are collected in Section~\ref{sec:proof}, with remaining proofs and supporting lemmas deferred to Appendix~\ref{sec:append}.

\subsection{Preliminaries}
We now introduce the preliminary background and notations. Let $d$ denote the dimension.
In this paper, unless otherwise specified, we restrict the domain to be $\Omega = [0,1]^d$ and the base measure $\pi$ to be the Lebesgue measure on $\Omega$. 
We use $\mu, \nu, \pi$ to denote probability distributions (measure), and also reserve $p_\mu(x), p_\nu(x), p_\pi(x)$ for the corresponding density functions w.r.t the Lebesgue measure (the Radon-Nikodym derivative). In other words, for ease of notation we use $\int_{\Omega} f(x) p_\mu(x) dx = \int_{\Omega} f d\mu$ to denote the same integration. $\| f\|_q := \left(\int_{\Omega} |f(x)|^q dx \right)^{1/q}$ denotes the $\ell_q$-norm, for $1\leq q \leq \infty$, and $\|w \|_q$ denotes the vector $\ell_q$-norm for a vector $\omega$. $[K]:=\{1, \ldots, K\}$ refers to the index set, for any $K \in \mathbb{N}_{>0}$. For a vector or a multi-index (possibly infinite-dimensional), the subscript $i$ denotes the $i$-th component.
We use the asymptotic notation $A(n) \precsim n^{\alpha}$, if $\limsup\limits_{n\rightarrow \infty} \log A(n)/\log n \leq \alpha$, holding other parameters fixed, similarly $A(n) \succsim n^{\alpha}$ if $\liminf\limits_{n\rightarrow \infty} \log A(n)/\log n \geq \alpha$. We call $A(n) \asymp n^{\alpha}$ when $A(n) \precsim n^{\alpha}$ and $A(n) \succsim n^{\alpha}$. 

Next, we introduce the function spaces. For a multi-index $\gamma \in \mathbb{N}_{\geq 0}^d$, we use $D^{(\gamma)} f$ to denote the $\gamma$-weak derivative of the function $f:\Omega \rightarrow \mathbb{R}$. For example, for the special case of infinitely smooth $f \in C^{\infty}(\Omega)$, $D^{(\gamma)} f$ takes the simple form $D^{(\gamma)} f = \partial^{|\gamma|} f/\partial x_1^{\gamma_1} \ldots \partial x_d^{\gamma_d}$, where $|\gamma| = \sum_i \gamma_i$.
\begin{definition}[Sobolev class: $\alpha \in \mathbb{N}_{>0}$]
	\label{def:sobolev}
	Given a smoothness parameter $\alpha \in \mathbb{N}_{\geq 0}$, $1\leq q \leq \infty$, and a radius $r \in \mathbb{R}_{\geq 0}$, the Sobolev class $W^{\alpha, q}(r)$ is defined as
	\begin{align*}
		W^{\alpha, q}(r) &:= \left\{ f \in \Omega \rightarrow \mathbb{R} : ( \sum_{|\gamma| \leq \alpha} \| D^{(\gamma)} f \|_{q}^q )^{1/q} \leq r \right\}, 
	\end{align*}
	where $\gamma$ is a multi-index and $D^{(\gamma)}$ denotes the $\gamma$-weak derivative. For the case $q=2$, we abbreviate the $W^{\alpha, q=2}(r)$ as $W^{\alpha}(r)$.
\end{definition}

To extend the analysis to distributions supported on a manifold $\Omega$, we further consider general Reproducing Kernel Hilbert Spaces (RKHS) $\cH \subset L^2_{\pi}$ (with $\pi$ as the base measure) endowed with the RKHS norm $\| \cdot \|_{\cH}$ and the corresponding positive semidefinite kernel $K(\cdot, \cdot): \Omega \times \Omega \rightarrow \mathbb{R}$. By the Mercer's theorem, one can characterize this RKHS via the following integral operator $\cT_\pi: L^2_\pi \rightarrow \cH$.

\begin{definition}[Integral operator of RKHS]
	\label{def:integral-operator}
	Define the integral operator $\cT_\pi: L^2_\pi \rightarrow \cH$,
	\begin{align*}
		\cT_{\pi} f(z) = \int_{\Omega} K(z, \cdot) f(\cdot) d \pi(\cdot), 
	\end{align*}
	and denote the eigenfunctions of this operator by $\psi_i$ and the associated eigenvalues by $t_i, i \in \mathbb{N}_{\geq 0}$ (sorted in a non-increasing order), with
	\begin{align*}
		\cT_{\pi} \psi_i = t_i \psi_i, ~\text{and}~ \int_{\Omega} \psi_i \psi_j d \pi = \delta_{ij}. 
	\end{align*}
\end{definition}
We assume that for all target distributions $\nu \in \G$ of interest and all $i\in \mathbb{N}_{\geq 0}$, there exists a universal constant on the variance of eigenfunctions in Def.~\ref{def:integral-operator}, 
\begin{align}
	\label{eq:assumption-rkhs}
	\sup_{\nu \in \G} \sup_{i\in \mathbb{N}_{\geq 0}} ~\E_{X \sim \nu} \psi_i(X)^2 \leq C.
\end{align}

To measure the complexity of functions from a learning theory perspective, we employ the following notion of combinatorial dimension for real-valued function, which is credited to \citet*{pollard1990empirical}. We will employ this combinatorial dimension as a complexity measure in deriving rates for neural network classes.
\begin{definition}[Pseudo-dimension]
	\label{def:comb-dim}
	Let $\F = \{ f: \Omega \rightarrow \mathbb{R} \}$ be a class of functions. The pseudo-dimension of $\F$, denoted by ${\rm Pdim}(\F)$, is the largest integer $m$ such that:
$\exists (X_i, y_i) \in \Omega \times \mathbb{R}, i\in [m]$, for any $(b_1, \ldots, b_m) \in \{ -1, 1 \}^m$ there exists $f \in \F$ such that $\sign(f(X_i) - y_i) = b_i, \forall i\in[m]$. 
\end{definition}

Finally, for two functions $f: \mathbb{R}^d \rightarrow \mathbb{R}$ and $g: \mathbb{R}^p \rightarrow \mathbb{R}^d$, denote $f\circ g$ to be the composition $f(g(x))$. We use the following notation for the composition of function classes
\begin{align}
\label{eq:composition}	
\F \circ \G := \{ f\circ g ~|~  f \in \F,  g \in \G\}.
\end{align}


\section{The Adversarial Framework}
\label{sec:nonparam}

We start by investigating the adversarial framework, including the Wasserstein, Sobolev, and MMD GAN, and GMM. Recall that the adversarial framework employed by GANs proposes to evaluate the accuracy of learning densities via the adversarial loss specified by the discriminator class. 
The goal of this section is to study the optimal minimax rates for learning a wide range of distributions, on a host of evaluation metric defined by the adversarial framework. Through the lens of nonparametric statistics, we answer how the structure of the distribution and the choice of the evaluation metric affects the optimal rates, and when fast rates are possible.

\subsection{Minimax Optimal Rates}
\begin{theorem}[Minimax optimal rates, Sobolev]
	\label{thm:optimal-rates-sobolev}
	Let $\Omega = [0,1]^d$.
	Consider the target distribution class $\G = \{ \nu~|~ p_\nu(x) \in  W^{\alpha}(r)\}$ in a Sobolev class with smoothness $\alpha \in \mathbb{N}_{\geq 0}$, for some constant $r>0$. Consider the evaluation metric induced by $\F = W^{\beta}(1)$, a Sobolev class with smoothness $\beta \in \mathbb{N}_{\geq 0}$. The minimax optimal rate is
	\begin{align*}
		\inf_{\widetilde{\nu}_n} \sup_{\nu \in \G} \E d_{\F}\left( \nu, \widetilde{\nu}_n \right) \asymp n^{-\frac{\alpha+\beta}{2\alpha+d}} \vee n^{-\frac{1}{2}},
	\end{align*}
	where $\widetilde{\nu}_n$ is any estimator for $\nu$ based on $n$ i.i.d. samples $X_1, X_2,\ldots X_n \sim \nu$.
\end{theorem}

\begin{remark}
	\rm
	
	The above establishes the minimax optimal rate for Sobolev GAN, with explicit dependence on the smoothness of the density $\alpha$ and that of the evaluation metric $\beta$ (here $\beta=1$ for Wasserstein GAN and $\beta = 0$ for TV as special cases, since our minimax lower bound holds for the special subclass $W^{\beta, \infty}$). First, note there is an interesting transition at $\beta = d/2$ (without depending on $\alpha$): above it the rate is parametric $n^{-\frac{1}{2}}$, and below it the rate is nonparametric. Second, to avoid the curse of dimensionality in the rates, one needs the sum of smoothness to be proportional to the dimension, that is, $\alpha+\beta = \Theta(d)$.  Note when $\beta$ is large, the rate is indeed faster, however, under a weaker evaluation metric. How to choose a good discriminator $\F$ in GANs with a provable guarantee under strong evaluation metrics such as $d_{TV}$ will be answered in Theorems \ref{thm:param-gans}-\ref{thm:leaky-ReLU}. 
	
\end{remark}

\begin{remark}[Relations to the literature]
	\rm
	The above theorem is an improvement to an earlier draft \citep*{liang2017well} of this paper, which was the first to formalize nonparametric estimation under the adversarial framework. The improvement on the upper bound is in one line of the original argument, specifically Eqn.~\eqref{eq:improvement}. The minimax lower bound of $n^{-\frac{\alpha+\beta}{2\alpha+d}}$ was first established in the earlier draft of this paper \cite[pages 18-19]{liang2017well}. In this version, we also provide a formal construction for the parametric lower bound of $n^{-\frac{1}{2}}$.
	We remark that the improvement of the upper bound in \cite*{liang2017well} was first derived in a follow-up work \citep{singh2018nonparametric} (see the discussion therein), in a general setting. Optimal upper bound of similar flavor was also obtained in \cite{mair1996statistical} with a different setup. After this paper was posted, there has been a growing list of works studying distribution estimation under the adversarial framework, with more general metrics and target distribution classes, to name a few, \cite{singh2018minimax,bai2018approximability,weed2019estimation,lei2019sgd, chen2020statistical}. 
	
	A closely related problem is: given $X_1, \ldots, X_n$ i.i.d. samples from $\nu$ and $Y_1,\ldots, Y_n$ i.i.d. samples from $\mu$, the optimal minimax rate for ``estimating the IPM'' between $\mu$ and $\nu$. In a follow-up paper, \cite{liang2019EstimatingCertain} showed that curiously, estimating the IPM itself is just as hard as estimating distributions under the IPM when $\beta<d/2$, in the following sense
	\begin{align}
		\frac{\log\log n}{\log n} \cdot n^{-\frac{\alpha+\beta}{2\alpha+d}} \precsim \inf_{\widetilde T_n} \sup_{\mu, \nu \in \G} \E |\widetilde T_n - d_{\F} (\mu, \nu) | \precsim n^{-\frac{\alpha+\beta}{2\alpha+d}} \enspace,
	\end{align}
	where $\widetilde{T}_n$ is any estimator based on the samples. This result shows that in the hard regime $\beta<d/2$, even evaluating whether two distributions are close is just as hard as estimating the distributions under the IPM $d_{\F}$.
	
\end{remark}

One can generalize the above theorem to more general RKHS. The motivation is to accommodate target distributions supported on image manifolds, with similarity better measured by non-linear kernels. It is useful to derive the explicit dependence on the intrinsic dimension of the manifold and the kernel, rather than the ambient dimension $d$. The Sobolev class considered in Thm.~\ref{thm:optimal-rates-sobolev} can be viewed as a special RKHS when the smoothness index is large enough \citep{devito2019ReproducingKernel}. In addition, the generalization will enable us to provide theoretical rates for MMD GAN \citep{dziugaite2015training,li2015generative,arbel2018gradient}.

\begin{theorem}[MMD rates, RKHS]
	\label{thm:optimal-rates-rkhs}
	Consider a RKHS $\cH \subset L^2_\pi$, with $\pi$ being a base measure. Assume that the eigenvalues of the integral operator $\cT_{\pi}$ decay as $t_i \asymp i^{-\kappa}$ for all $i \in \mathbb{N}_{\geq 0}$, with parameter $\kappa \in \mathbb{R}_{>0}$. Consider the evaluation metric induced by $\F = \{ f ~|~ \| f \|_{\cH} \leq 1\}$, and the target distribution $\nu$ whose Radon-Nikodym derivative $\frac{d\nu}{d\pi}$ w.r.t $\pi$ lies in a smooth subclass $\G = \{ \nu ~|~ \| \cT_{\pi}^{-(\lambda-1)/2} \frac{d\nu}{d\pi} \|_{\cH} \leq r  \}$ with smoothness parameter $\lambda \in \mathbb{R}_{>0}$ (for some fixed radius $r>0$). Under the assumption~\eqref{eq:assumption-rkhs}, the following upper rate holds
	\begin{align*}
		\sup_{\nu\in \G} \E d_{\F}\left( \nu, \widetilde{\nu}_n \right) \precsim n^{-\frac{(\lambda+1)\kappa}{2\lambda \kappa+2}} \vee n^{-\frac{1}{2}} .
	\end{align*}
\end{theorem}
\begin{remark}[Intrinsic dimension]
	\rm
	Remark that the above corollary works with general base measure $\pi$ and domain $\Omega$.
	Here the target (Radon-Nikodym derivative $\frac{d\nu}{d\pi}$) lies in a subclass of the RKHS when $\lambda>1$, with $\lambda$ quantifying its smoothness: the high frequency component decays sufficiently fast. This is a standard formulation studied in the RKHS literature, see \cite*{caponnetto2007optimal}. The parameter $\kappa$ describes the intrinsic dimension of the integral operator. When $\kappa > 1$, the intrinsic dimension (trace of $\cT_\pi$) is bounded as ${\rm Tr}(\cT_\pi) =  \sum_{i \geq 1} i^{-\kappa} \leq C$, therefore the upper bound reads the parametric rate $n^{-\frac{(\lambda+1)\kappa}{2\lambda \kappa+2}} \vee n^{-\frac{1}{2}} = n^{-\frac{1}{2}}$. When $\kappa < 1$, to obtain 
	$
	\E d_{\F}\left( \nu, \widetilde{\nu}_n \right) \leq \epsilon,
	$
	the sample complexity scales as $$n \asymp \epsilon^{2 + \frac{2}{\lambda+1}\left( \frac{1}{\kappa} - 1 \right) }.$$ Therefore the curse of dimensionality only appears in the effective dimension, described by $1/\kappa$.

	The Sobolev class $W^{\beta}$ can be regarded as a special RKHS with $\kappa = 2\beta/d$. The reason can be seen from the Sobolev ellipsoid Def.~\ref{eq:sobolev-ellipsoid} (a weighted $L^2$-space same as RKHS), with corresponding a weight-decay $t_i \asymp (i^{2/d})^{-\beta} \asymp (1+\| \xi \|^2)^{-\beta}$. Here $i(\xi)$ is the lexicographic re-ordering of the multivariate Fourier index $\xi$. In such a case, the generator class $W^{\alpha}$ can be thought as subclass $\G$ in Thm.~\ref{thm:optimal-rates-rkhs} with $\lambda = \alpha/\beta$. Plug in $\lambda = \alpha/\beta$ and $\kappa = 2\beta/d$, the bound in Thm.~\ref{thm:optimal-rates-rkhs} recovers $n^{-\frac{\alpha+\beta}{2\alpha+d}}$ which agrees with Thm.~\ref{thm:optimal-rates-sobolev}. Therefore the lower bound in Thm.~\ref{thm:optimal-rates-sobolev} suggests that the rate for MMD GAN is also sharp, for a particular subclass. 
\end{remark}

\subsection{Oracle Inequality and Regularization}
\label{sec:nonparam-oracle}

In this section, we use a simple oracle inequality to show that when the generator class $\D_G$---typically represented by neural networks---is misspecified for the target distribution $\nu$, one can still derive oracle results based on the adversarial framework. 

Let us recall the notations. Denote $\D_{G}$ to be class of distributions represented by the generator, and $\F_{D}$ to be the class of functions realized by the discriminator
\begin{align}
	\label{eq:gan-nonparam}
	\mu_n = \argmin_{\mu \in \D_{G}} \max_{f \in \F_{D}} \left\{ \E_{Y \sim \mu} f(Y) - \E_{X \sim \nu_n} f(X) \right\}.
\end{align}
where $\nu_n$ is some estimate of the distribution based on $n$ i.i.d. drawn samples from the target distribution $\nu$.

The goal in this section to extend our adversarial framework to obtain upper rates for \eqref{eq:gan-nonparam}. In addition, the oracle inequalities (Lemma \ref{lem:oracle-ineq} and \ref{lem:oracle-ineq-general}) developed will be crucial for model misspecification, which makes the results of practical relevance.

\begin{theorem}[Misspecification: nonparametric]
	\label{thm:nonparam-gans}
	Let $\D_G$ be any generator class. 
	Consider the discriminator metric induced by $\F_{D} = W^{\beta}(1)$. Consider the target density $p_\nu \in  W^{\alpha}(r)$.
	With the empirical distribution $\widehat{\nu}^{n} := \frac{1}{n} \sum_{i=1}^n \delta_{X_i}$ as the plug-in, the GAN estimator 
	\begin{align*}
		\widehat{\mu}_{n} \in \argmin_{\mu \in \D_G} \max_{f \in \F_D} \left\{ \int_{\Omega}  f d\mu    -  \int_{\Omega} f d\widehat{\nu}^{n}  \right\} \;,
	\end{align*}
	learns the target distribution with rate
	\begin{align*}
		\E d_{\F_D}(\widehat{\mu}_n, \nu) \leq  \min_{\mu \in \D_G} d_{\F_D}(\mu, \nu) +   n^{-\frac{\beta}{d}} \vee \frac{\log n}{\sqrt{n}} \;.
	\end{align*}
	In contrast, there exists a regularized $\widetilde{\nu}^{n}$ as the plug-in
	\begin{align*}
		\widetilde{\mu}_{n} \in \argmin_{\mu \in \D_G} \max_{f \in \F_D}  \left\{ \int_{\Omega}  f d \mu   -  \int_{\Omega} f d\widetilde{\nu}^{n}  \right\} \;,
	\end{align*}
	where a faster rate is attainable
	\begin{align*}
		\E d_{\F_D}(\widetilde{\mu}_n, \nu) \leq  \min_{\mu \in \D_G} d_{\F_D}(\mu, \nu) +   n^{-\frac{\alpha+\beta}{2\alpha+d}} \vee \frac{1}{\sqrt{n}} \;.
	\end{align*}
\end{theorem}

The proof of the above theorem is based on a simple oracle inequality Lemma~\ref{lem:oracle-ineq}. Later, we will generalize the oracle inequality (see Lemma~\ref{lem:oracle-ineq-general}) to establish rates when both the generator and discriminator are neural networks, and when one only has finite $m$ simulated samples from the generator. Curiously, a generalization of the oracle inequality gives rise to a curious notion of pair regularization, which we will study in Section~\ref{sec:param}.

\begin{remark}[Regularization]
	\rm
Observe that the rates satisfy $n^{-\frac{\alpha+\beta}{2\alpha+d}} \vee n^{-1/2} \precsim n^{-\frac{\beta}{d}} \vee n^{-1/2}\log n$. Namely, the regularized empirical distribution as the plug-in for GANs attains a better upper bound. In a high level, the regularized empirical distribution $\widetilde \nu^n$ filters out high frequency component of the empirical distribution to enforce regularization.
We mention that to obtain an implementable algorithm for the smoothed/regularized empirical distribution $\widetilde{\nu}^n$ in Thm.~\ref{thm:nonparam-gans}, one may use the following in practice
\begin{align*}
	\frac{d\widetilde{\nu}^n}{dx} = \frac{1}{nh_n}\sum_{i\in [n]} K\left(\frac{x-X_i}{h_n}\right),
\end{align*}
with specific choices of the kernel $K$ and bandwidth $h_n$. When using the Gaussian kernel, this so-called ``instance noise'' technique \citep{sonderby2016amortised,arjovsky2017towards,mescheder2018training} is used in GAN training: each time when evaluating the stochastic gradients for generator and discriminator, sample a mini-batch of data and then perturb them by a Gaussian. Statistically, one may view this data augmentation (or stability to data perturbation) as a form of regularization \citep*{yu2013stability}, to prevent the generator from memorizing the empirical data and learning a too complex model. We will show in Section~\ref{sec:param} that, curiously, the choice of generator and discriminator pair can also serve the goal of regularization. 
\end{remark}

\subsection{Application: Generalized/Simulated Method of Moments}

We show how the adversarial framework and the oracle inequality established so far can be used to derive non-asymptotic bounds for using the Generalized Method of Moments (GMM) to estimate distributions implicitly \citep{hansen1982LargeSample,back1993ImpliedProbabilities}, with a particular choice of $\F_{D}$ specifying the moment conditions.

Let $\{ \phi_k(x), k\in [K]\}$ be a set of functions of cardinality $K$ that determine the moment conditions. For instance, $\phi_k(x) = x^k, k\in [K]$ in the standard method of moments. Let $\bW \in \mathbb{R}^{K \times K}$ be a symmetric positive-definite matrix. Specify the discriminator class as
\begin{align}
	\label{eq:GMM-discriminator}
	\F_D = \left\{ f(x) = \sum_{k \in [K]} \omega_k \phi_k(x) ~|~ \omega^\top \bW^{-1} \omega \leq 1 \right\}.
\end{align}
In such a case, the adversarial framework reduces to GMM.
\begin{corollary}[GMM and the Adversarial Framework]
	\label{coro:GMM}
	Let $K\in \mathbb{N}$ be finite and $\bW \in \mathbb{R}^{K\times K}$ be a symmetric positive-definite weight matrix.  With the choice of $\F_D$ in \eqref{eq:GMM-discriminator}, the GAN estimator $\widehat\mu_n$ in \eqref{eq:gan-nonparam}
	implements the GMM. In addition, the following non-asymptotic guarantee holds
	\begin{align*}
		\E d_{\F_D}(\widehat{\mu}_n, \nu) \leq  \min_{\mu \in \D_G} d_{\F_D}(\mu, \nu) + 4 \cdot \sqrt{\frac{\E_{X\sim \nu}[\sum_{i,j\in [K]} \bW_{ij} \phi_i(X) \phi_j(X) ]}{n}} \enspace.
	\end{align*}
\end{corollary}
Remark that under mild conditions on the moments $\E_{X\sim \nu} \phi_k(X)^2, k\in [K]$ and the weight matrix $\bW$, the non-asymptotic bound is of the order $\sqrt{K/n}$, which is the optimal parametric rate $\sqrt{\text{Pdim}(\F_{D})/n}$, in view of the lower bound in Thm.~\ref{thm:optimal-rates-sobolev}.

Let us elaborate on why the Simulated Method of Moments (SMM) is contained in the GANs formulation to be considered in Section~\ref{sec:param}. Recall the GANs estimator \eqref{eq:gan-network} with parameter $\theta$, the generalized moment equations $M(\theta, x) \in \mathbb{R}^K$ in such case are precisely
\begin{align}
	M_k(\theta, x) = \E_{Z\sim \pi} \left[\phi_k(g_{\theta}(Z))\right] - \phi_k(x), ~k\in [K]
\end{align}
with the simulated counterpart
\begin{align}
	M_k^{\text{sim}}(\theta, x) =  \widehat\E_{m} \left[\phi_k(g_{\theta}(Z))\right] - \phi_k(x) \enspace.
\end{align}
The SMM is the GANs estimator as in \eqref{eq:gan-network}
\begin{align}
	\widehat{\theta}_{m,n} \in \argmin_{\theta: g_\theta \in \G} \left(\widehat\E_n M^{\text{sim}}(\theta, X)\right)^\top \bW \left(\widehat\E_n M^{\text{sim}}(\theta, X)\right) \enspace.
\end{align}
Thm.~\ref{thm:param-gans} and \ref{thm:param-gans-hellinger} in the next section derive non-asymptotic upper bounds for such SMM, with the form $\sqrt{\text{Pdim}(\F)/n} + \sqrt{\text{Pdim}(\F\circ \G)/m}$ (overlooking logarithmic factors).


\section{Generative Adversarial Networks}
\label{sec:param}

In this section, we consider when both the generator and discriminator are parameterized by neural networks, and derive rates applicable to GAN used in practice (as well as SMM). To be specific, let $\F = \{ f_\omega(x): \mathbb{R}^d \rightarrow \mathbb{R} \}$ be the discriminator functions realized by a neural network with parameter $\omega$ describing the weights of the network. Let $\G = \{ g_{\theta}(z): \mathbb{R}^{d} \rightarrow \mathbb{R}^{d} \}$ be the generator neural network transformation with weights parameter $\theta$. We keep this parametrization in an abstract form for now as we will first state our general theorems before applying them to specific cases. Consider $Z \sim \pi$ as the random input with distribution $\pi$, and the target distribution $X \sim \nu$. Denote $\mu_\theta$ as the probability distribution of $g_{\theta}(Z)$. Consider the parametrized GAN estimator used in practice
	\begin{align}
		\label{eq:gan-nn}
		\widehat{\theta}_{m,n} \in \argmin_{\theta: g_\theta \in \G} \max_{\omega: f_\omega \in \F} ~~ \left\{ \widehat{\E}_m  f_\omega (g_\theta(Z))  - \widehat{\E}_n  f_\omega(X)  \right\},
	\end{align}
where $m$ and $n$ denote the number of the generator samples (simulated) and target distribution samples.

Let us state the goal of the current section and connect to the adversarial framework established. So far, we have derived the optimal rates for nonparametric densities under strong evaluation metrics such as Wasserstein ($\beta=1$) or total variation distance ($\beta=0$). The curse of dimensionality in the sample complexity is inevitable unless the distribution class of interest is sufficiently structured (smooth). Two questions naturally arise. First, for the structured parametric distributions such as those parametrized by the generator networks in GANs, are fast parametric rates attainable? Second, can one obtain fast rates under the strong evaluation metric using a discriminator metric induced by neural networks as in GANs, which misspecifies the evaluation metric? We will answer both questions, directly for GANs estimator \eqref{eq:gan-nn}.

\subsection{Generalized Oracle Inequality and Parametric Rate}

First, we will generalize the oracle inequality to GANs estimator $\widehat{\theta}_{m,n}$ \eqref{eq:gan-nn}. Then we will show that the oracle approach, when applied to neural network classes, sheds light on the choice of \textit{generator-discriminator-pair} as regularization. 

\begin{lemma}[Generalized oracle inequality]
	\label{lem:oracle-ineq-general}
	Consider the GANs estimator $\widehat{\theta}_{m,n}$ defined in \eqref{eq:gan-nn}.
	Recall the composition in Def. \eqref{eq:composition}. For any $g_\theta \in \G$, denote $\mu_\theta$ as the probability distribution of $g_{\theta}(Z), Z\sim \pi$.
	Under the condition that $\F$ and $\F \circ \G$ are symmetric,
	the following oracle inequality holds for any $\theta$ with $g_\theta \in \G$,
	\begin{align*}
		d_{\F}\left( \mu_{\widehat{\theta}_{m,n}}, \nu  \right) \leq d_{\F}(\mu_\theta, \nu) + 2 d_{\F}\left( \widehat{\nu}^n, \nu \right) + d_{\F}\left( \widehat{\mu}^m_\theta, \mu_\theta \right) + d_{\F \circ \G}(\widehat{\pi}^m, \pi).
	\end{align*}
	Here for any distribution $\mu$, $\widehat{\mu}^n$ denotes the empirical distribution based on $n$ i.i.d. samples from $\mu$.
\end{lemma}

The innovative aspects of the above Lemma are two-fold. Firstly, the Lemma provides upper bound on the 
\textit{implicit} distribution estimator $\mu_{\widehat{\theta}_{m,n}}$ (distribution of the random variable $g_{\widehat{\theta}_{m,n}}(Z)$), without knowing the explicit form of the density function in general. Note that we do have direct sampling mechanisms by transforming the random variable $Z$, which is a computational advantage. Secondly, we make explicit the dependence on the number of generator samples $m$, in addition to the number of target samples $n$. The role and complexity of the generator network is made explicit in the bound.  It is clear that when $m \rightarrow \infty$, the current lemma reduces to Lemma~\ref{lem:oracle-ineq}. This Lemma made explicit the choice of simulated samples $m$ relative to the real-data samples $n$, and how the classes $\G$ and $\F$ affect the equation. Clearly, such a bound will also prove useful in analyzing SMM.

Next, we apply Lemma~\ref{lem:oracle-ineq-general} to establish parametric rates for distributions realized by neural networks, in the following Thm.~\ref{thm:param-gans} and \ref{thm:leaky-ReLU} (with their corollaries). 
We emphasize again here that GANs only use a misspecified discriminator $\F$ parametrized by neural networks with limited capacity. And $d_{\F}$ is \textit{different} from the the objective evaluation metrics such as $d_{TV}, d_{H}$.

\begin{theorem}[GANs upper rate on KL: parametric]
	\label{thm:param-gans}
	Consider GANs estimator
	\begin{align}
		\label{eq:gan-nn-B}
		\widehat{\theta}_{m,n} \in \argmin_{\theta: g_\theta \in \G} \max_{\omega: f_\omega \in \F, \| f_\omega \|_\infty \leq B} ~~ \left\{ \widehat{\E}_m  f_\omega (g_\theta(Z))  - \widehat{\E}_n  f_\omega(X)  \right\} \;.
	\end{align}
	where $B>0$ is some absolute constant, $m$ and $n$ denote the number of the generator samples and target distribution samples. Recall the pseudo-dimension defined in Def.~\ref{def:comb-dim}. For total variation distance, and Kullback-Leibler divergence, we have 
	\begin{align*}
		&\E d_{TV}^2\left(\nu, \mu_{\widehat{\theta}_{m,n}} \right)  \leq \frac{1}{4}\left[ \E d_{KL}\left(\nu ||  \mu_{\widehat{\theta}_{m,n}} \right) + \E d_{KL}\left(\mu_{\widehat{\theta}_{m,n}} || \nu \right) \right]   \\
		& \leq \frac{1}{2} \cdot \sup_{\theta} \inf_{\omega} \left\| \log \frac{p_\nu}{p_{\mu_{\theta}}} - f_\omega \right\|_{\infty} + \frac{B}{4\sqrt{2}}\cdot  \inf_\theta  \left\| \log \frac{p_{\mu_\theta}}{p_{\nu}} \right\|_{\infty}^{1/2}  \\
		& \quad \quad +  C \cdot  \sqrt{\text{Pdim}(\F) \left( \frac{\log m}{m} \vee \frac{\log n}{n} \right)} \vee  \sqrt{  \text{Pdim}(\F \circ \G) \frac{\log m}{m} } \; , 
	\end{align*} 
	where $C>0$ is some universal constant independent of $\text{Pdim}(\F)$, $\text{Pdim}(\F \circ \G)$ and $m, n$. 
\end{theorem}
Note that when $\F, \G$ are both neural network classes, $\F\circ \G$ is also a neural network class with a larger architecture.
The upper bound above on the Jensen-Shannon, Kullback-Leibler divergence, and TV distance consists of three parts: the approximation errors $A_1(\F,\G,\nu)$, $A_2(\G, \nu)$ and the stochastic error $S(\F, \G, n, m)$,
\begin{align}
	A_1(\F,\G,\nu) &:= \frac{1}{2} \cdot \sup_{\theta} \inf_{\omega} \left\| \log \frac{p_\nu}{p_{\mu_{\theta}}} - f_\omega \right\|_{\infty} \;, \nonumber \\
	A_2(\G, \nu) & := \frac{B}{4\sqrt{2}} \cdot  \inf_\theta  \left\| \log \frac{p_{\mu_\theta}}{p_{\nu}} \right\|_{\infty}^{1/2} \;, \label{eq:approx-stat-error}  \\
	S_{n,m}(\F, \G) &:=  \sqrt{\text{Pdim}(\F) \left( \frac{\log m}{m} \vee \frac{\log n}{n} \right)} \vee  \sqrt{  \text{Pdim}(\F \circ \G) \frac{\log m}{m} } \nonumber \;. 
\end{align}
We emphasize that the term $A_1(\F,\G,\nu)$ is in a maximin form $\sup_{\theta} \inf_{\omega}$, which is crucial and differs from the minimax form $\inf_{\theta} \sup_{\omega}$. In English, $A_1(\F,\G,\nu)$ describes how the best discriminator function $f_\omega$ that can express the class of density ratios $p_{\mu_{\theta}}/p_{\nu}$, $A_2(\G, \nu)$ reflects the the expressiveness of the generator class, and $S_{n,m}(\F, \G)$ describes the statistical complexity of both the generator and discriminator.
In the next section, we will elaborate on the interplay among the two approximation error terms $A_1(\F,\G,\nu), A_2(\G, \nu)$, and the stochastic error term $S_{n,m}(\F, \G)$.
\begin{remark}
	\rm
	To obtain non-trivial rates, the above theorem requires $\mu_\theta$ and $\nu$ to be absolutely continuous, for all $\theta$ of interest. However, this is not essential, as similar results hold qualitatively the same for the non-absolutely continuous case, based on the Hellinger distance. As shown in the next theorem, 
	$-1\leq \frac{\sqrt{p_\nu} - \sqrt{p_{\mu_{\theta}}}}{\sqrt{p_\nu} + \sqrt{p_{\mu_{\theta}}}} \leq 1$ is well-defined even for non-absolutely continuous distributions $\mu_\theta$ and $\nu$. 
\end{remark}

\begin{theorem}[GANs upper rate on Hellinger: parametric]
	\label{thm:param-gans-hellinger}
	Consider the same GANs estimator $\widehat{\theta}_{m,n}$ as in Thm.~\ref{thm:param-gans}. The for the Hellinger distance, 
	\begin{align*}
		d_{H}(\mu, \nu) := \left( \int \left(\sqrt{p_\mu} - \sqrt{p_\nu} \right)^2 dx \right)^{1/2},
	\end{align*}
	we have
	\begin{align*}
			&\E d_{TV}^2\left(\nu, \mu_{\widehat{\theta}_{m,n}} \right)  \leq \E d_{H}^2 \left(\nu, \mu_{\widehat{\theta}_{m,n}} \right)   \\
			& \leq 2\cdot \sup_{\theta} \inf_{\omega} \left\| \frac{\sqrt{p_\nu} - \sqrt{p_{\mu_{\theta}}}}{\sqrt{p_\nu} + \sqrt{p_{\mu_{\theta}}}} - f_\omega \right\|_{\infty} + 2B \cdot  \inf_\theta  \left\| \frac{\sqrt{p_\nu} - \sqrt{p_{\mu_{\theta}}}}{\sqrt{p_\nu} + \sqrt{p_{\mu_{\theta}}}} \right\|_{\infty} \\
			& \quad \quad +  C \cdot  \sqrt{\text{Pdim}(\F) \left( \frac{\log m}{m} \vee \frac{\log n}{n} \right)} \vee  \sqrt{  \text{Pdim}(\F \circ \G) \frac{\log m}{m} } , 
	\end{align*} 
		where $C>0$ is some universal constant.
\end{theorem}

Finally, as a corollary of Thm.~\ref{thm:param-gans}, one can establish similar results for the Wasserstein distance. 
\begin{corollary}
	\label{cor:wass}
	Recall the definitions in \eqref{eq:approx-stat-error}. 
	Assume that $\F$ is with Lipschitz constant $L_\F$ and $\G$ with $L_\G$. Then for either
	(1) $Z \sim N(0, I_d)$, or (2) $Z, X$ lie in $[0, 1]^d$, we have
	\begin{align*}
		\E d_{W}^2\left(\nu, \mu_{\widehat{\theta}_{m,n}} \right) &\leq C_1 \cdot A_1(\F,\G,\nu)+ C_2 \cdot A_2(\G, \nu) + C_3 \cdot S_{n,m}(\F, \G)
	\end{align*}
	where $C_1, C_2, C_3>0$ are some constants independent of $\text{Pdim}(\F)$, $\text{Pdim}(\F \circ \G)$ and $m, n$, but depend on $L_\F, L_\G$.
\end{corollary}

\subsection{Generator-Discriminator-Pair Regularization}
\label{sec:pair-regularization}

\begin{figure}[ht!]
\centering
\includegraphics[width=0.6\textwidth]{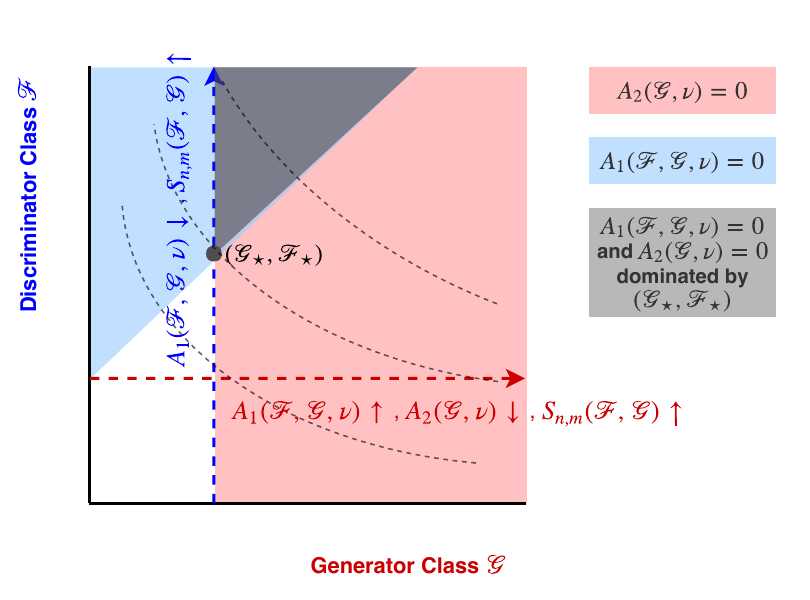}
\caption{Pair regularization diagram on how well GANs learn distributions in TV distance, when tuning with generator $\G$ and discriminator $\F$ pair. 
The diagram is illustrated based on upper bounds on TV distance, namely $A_1(\F,\G,\nu)+ A_2(\nu, \G) +S_{n,m}(\F, \G)$ in Thm.~\ref{thm:param-gans}.
The red shaded region corresponds to $A_2(\G, \nu) = 0$ and the blue shaded region is $A_1(\F,\G,\nu) = 0$. The grey dashed line corresponds to the indifference curve for the statistical error $S_{n,m}(\F, \G)$. One can see that the choice $(\G_\star, \F_\star)$ dominates the other choices in the grey shaded area, and the other choice on the same grey dashed line.}
\label{fig:pair-dis-gen}
\end{figure}

In this section, we discuss the new pair regularization, and its trade-offs presented in Thm.~\ref{thm:param-gans}.
In GANs, both the generator and discriminator are choices of tuning parameters for users to specify. Therefore, the trade-off between approximation error and stochastic error is more complicated. We use the following two thought experiments to explain the intricacies of the generator-discriminator-pair choice. 
\begin{enumerate}
	\item For a fixed generator class $\G$, when the discriminator class $\F$ increases the complexity, it will be easier for the discriminator to tell apart good and bad generators in the TV sense (w.r.t. the target distribution). However, the stochastic error becomes larger as one is learning from a large discriminator model in GANs. This is reflected in the upper bounds obtained in Thm.~\ref{thm:param-gans} and \ref{thm:param-gans-hellinger}, shown along the blue dashed arrow direction in Fig.~\ref{fig:pair-dis-gen}. 
	\item For a fixed discriminator class $\F$, as the generator $\G$ becomes more complex, it is capable of expressing distributions closer to the target. However, at the same time, it introduces difficulty for two reasons. First, the generator may create distributions that are far from the target distribution in the TV sense, but indistinguishable to the discriminator. Second, the stochastic error becomes larger as one is learning from a larger generator model. The above is shown by the red dashed arrow direction in Fig.~\ref{fig:pair-dis-gen}.
\end{enumerate}
In general, regularization using the generator-discriminator-pair is more subtle than the conventional approximation and stochastic error (or bias and variance) trade-off. We visualize such trade-offs in Fig.~\ref{fig:pair-dis-gen}, with $A_1(\F,\G,\nu)$, $A_2(\G, \nu)$ and $S_{n,m}(\F, \G)$ defined in \eqref{eq:approx-stat-error}. Here, the tuning parameters lie in a two-dimensional domain, rather than in a one-dimensional index. For a fixed target $\nu$, as $(\G, \F)$ both become richer, $A_2(\G, \mu)$ decreases, $S_{n,m}(\F, \G)$ increases, but $A_1(\F,\G,\nu)$ may increase, decrease or stay unchanged. On the one hand, one can eliminate some $(\G, \F)$ pairs using notions of dominance. The simple U-shaped picture for bias-variance trade-off no longer exists. On the other hand, by stepping into the two-dimensional tuning domain, there are more choices for tuning pairs that potentially give rise to better rates, which we will showcase in Thm.~\ref{thm:leaky-ReLU}.

The following corollary concerns $A_1(\F,\G,\nu)$ and $A_2(\G, \nu)$ through choosing the generator-discriminator-pair, as a step towards understanding the new notion of pair regularization for GANs.
\begin{corollary}[Choice of generator and discriminator]
	\label{cor:design}
	Consider the target distribution class $\nu \in \D_{R}$, and the generator distribution class $\mu_{\theta} \in \D_G$. With the discriminator chosen as
	\begin{align*}
		\F_D  := \{ \log(p_\nu) - \log(p_{\mu_{\theta}}) ~|~ \text{for all}~ \nu \in \D_R,~  \mu_{\theta} \in \D_G\},
	\end{align*}
	then
	\begin{align}
		\label{eq:choice-dis}
		A_1(\F,\G,\nu) = 0.
	\end{align}
	In addition, if the generator is well-specified in the sense $\D_G \supseteq \D_R$, then
	\begin{align}
		\label{eq:choice-gen}
		A_2(\G,\nu) = 0.
	\end{align}
	With such choice of $\D_G$ and $\F_D$,
	$\E d_{TV}^2\left(\nu, \mu_{\widehat{\theta}_{m,n}} \right)  \precsim S_{n,m}(\F, \G).$
\end{corollary}

\begin{remark}[Pair regularization diagram]
\rm

Let us illustrate the above corollary using Fig.~\ref{fig:pair-dis-gen}. Eqn.~\eqref{eq:choice-dis} corresponds to the blue shaded region in the diagram, Eqn.~\eqref{eq:choice-gen} represents the red shaded region, and the intersection is highlighted by the grey shaded region. At the intersection, the approximation error $A(\F,\G,\nu)$ is zero, so all pairs are dominated by the choice $(\G_\star, \F_\star)$ (as other pairs have a larger variance $S_{n,m}(\F, \G)$). In addition, we argue that $(\G_\star, \F_\star)$ is also the best solution along the indifference curve for $S_{n,m}(\F, \G)$, denoted by the grey dashed line. To see this, moving $(\G_\star, \F_\star)$ towards the northwest direction on the indifference curve away from $(\G_\star, \F_\star)$, $A_1, S_{m,n}$ stay unchanged, but $A_2(\G_\star, \nu) \leq  A_2(\G', \nu)$. Moving $(\G', \F')$ towards the southeast direction, $A_2, S_{m,n}$ stay the same, but $A_1(\G_\star, \F_\star, \nu) \leq A_1(\G', \F', \nu)$. Similarly, one can argue that all pairs above the indifference curve is dominated by $(\G_\star, \F_\star)$. 
We defer the further discussion on the pair regularization versus the classic regularization to Section~\ref{sec:discussion}. 

\end{remark}

\subsection{Applications: Deep ReLU Networks and More}
We showcase how to apply our pair regularization theory to GANs used in practice in this section. We consider two special cases of neural network generator and discriminator and derive the rates for implicitly estimating certain parametric distributions. The two applications demonstrate two extreme cases of GANs. One is a deep ReLU GAN for learning a complex implicit distribution. The other is a two-layer GAN for learning a simple multivariate Gaussian distribution. In both cases, our theory exhibits a near-optimal dependence on the dimension and the network depth.

Let's introduce the neural networks parameter space.
The \textit{generator} $g_{\theta}(\cdot) : z\mapsto x$ is parametrized by a Multi-Layer Perceptron (MLP):
	\begin{align*}
		h_0 &= z, \\
		h_l &= \sigma_a(W_l h_{l-1} + b_l), ~ 0< l <L \\
		x   &= W_L h_{L-1} + b_L,
	\end{align*}
	where $h_l$ denotes the output of hidden units, and $x$ is the final output of the MLP. Here the activation is leaky ReLU
\begin{align}
	\label{eq:leaky-ReLU}
	\sigma_a(t) = \max\{ t, a t \},~ \text{for some fixed $0<a\leq 1$}.
\end{align} 
	Denote the parameter space for the generator weights as
	\begin{align*}
		\theta \in \Theta(d, L) := \{ \theta= (W_l \in \mathbb{R}^{d \times d}, b_l \in \mathbb{R}^d, 1\leq l \leq L) ~|~ {\rm rank}(W_l) = d, \forall 1\leq l\leq L  \}.
	\end{align*}
We require the $W_l$ to be full rank so that the generator transformation $g_{\theta}$ is invertible. 
One can verify, when the input distribution $Z \sim U([0,1]^d)$ is uniform, the class of densities realizable by $g_{\theta}(Z)$, for $\theta \in \Theta(d, L)$ has the following closed form, 
\begin{align}
	\label{eq:analytic-form}
	\log \big(p_{\mu_\theta}(x) \big) = c_1 \sum_{l=1}^{L-1} \sum_{i=1}^d \mathbbm{1}_{m_{li}(x) \geq 0} + c_0(\theta),
\end{align}
with some proper choice of $c_1, c_0(\theta)$. Here $m_{li}(x)$ is the function computed by the $i$-th hidden unit in the $l$-th layer of a certain MLP \footnote{The architecture and weights depend on the generator network $g_\theta$, with depth $L$ and $d$ hidden units in each layer.}, with dual leaky ReLU activation (defined in next paragraph) and weights properly chosen as a function of $\theta$. For details, see derivation \eqref{eq:m} and \eqref{eq:realizable-density-MLP}. Remark that from the closed form expression, as the depth grows (as a function of $n$), the generator is capable of expressing increasingly complex distributions. Clearly from the expression, one can see that for any $\theta, \theta' \in \Theta(d, L)$, $\mu_\theta$ and $\mu_{\theta'}$ are absolutely continuous.

The \textit{discriminator} $f_\omega(x): \mathbb{R}^d \rightarrow \mathbb{R}$ is parametrized by a feedforward neural network with activation functions include dual leaky ReLU activation
\begin{align}
\label{eq:leaky-ReLU-dual}
\sigma^\star_a(t) := \min\{ t, a t \}, ~\text{for $a \geq 1$},
\end{align} 
and threshold activation $\sigma^\star_\infty(t) := \mathbbm{1}_{t\leq 0}$. 
The feedforward network has the following structure: hidden units are grouped in a sequence of $L$ layers (the depth of the network), where a node is in layer $1\leq l \leq L$, if it has a predecessor in layer $l-1$ and no predecessor in any layer $l' \geq l$. Computation of the final output unit proceeds layer-by-layer: at any layer $l < L$, each hidden unit $u$ receives an input in the form of a linear combination $\widetilde{x}_u' w_u + b_u$, and then outputs $\sigma_a(\widetilde{x}_u' w_u + b_u)$, where the vector $\widetilde{x}_u$ collects the output of all the units with a directed edge into $u$ (that is, from prior layers). $\omega$ denotes all the weights in such feedforward network.

\begin{figure}[ht!]
\centering
\includegraphics[width=0.6\textwidth]{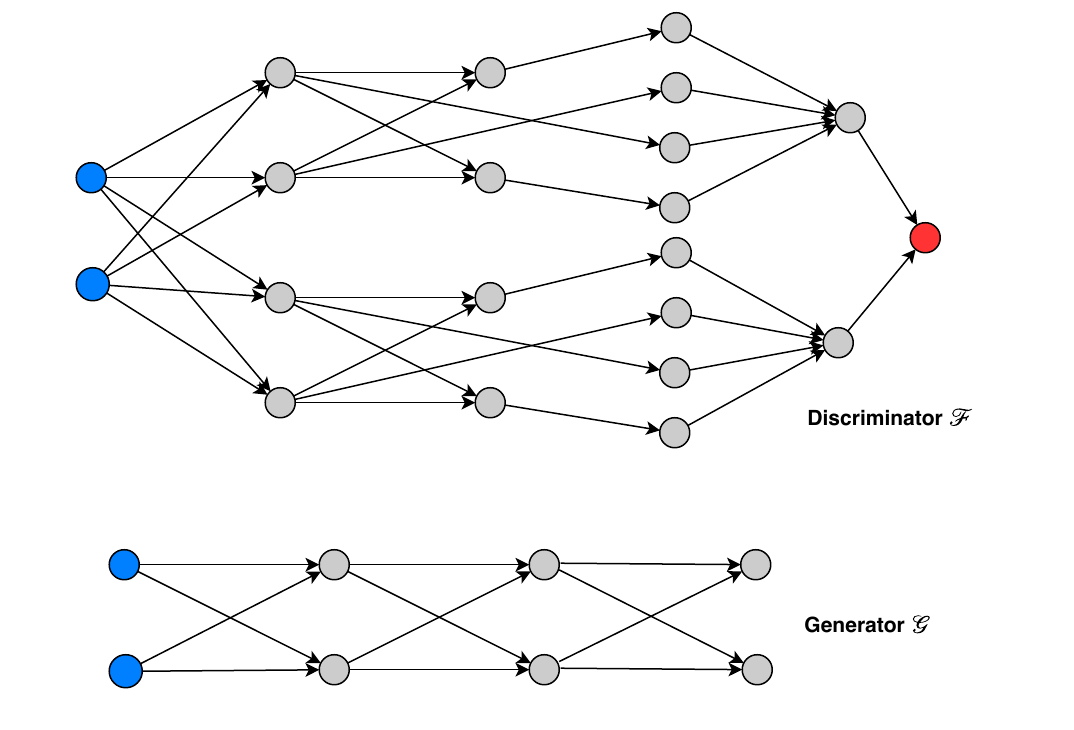}
\caption{Illustration of discriminator $\F$ (feed-forward network) and generator $\G$ (multi-layer perceptron) in Thm.~\ref{thm:leaky-ReLU}, for $L=3$. }
\label{fig:relu}
\end{figure}

\begin{theorem}[Leaky-ReLU generator and discriminator]
	\label{thm:leaky-ReLU}
	
	Consider a MLP generator $g_{\theta}:\mathbb{R}^d \rightarrow \mathbb{R}^d$, $\theta \in \Theta(d, L)$ with depth $L$ and width $d$, using leaky ReLU $\sigma_a(\cdot)$ activation \eqref{eq:leaky-ReLU} with any $0<a \leq 1$. Consider the class of realizable distributions, that is, $X \sim \nu$ enjoys the same distribution as $g_{\theta_*}(Z)$ with some $\theta_* \in \Theta(d, L)$ and $Z \sim U([0,1]^d)$. Choose the discriminator $f_\omega: \mathbb{R}^d \rightarrow \mathbb{R}$ to be a feedforward neural network (architecture shown in Fig.~\ref{fig:relu}) with depth $L+2$, using dual leaky ReLU $\sigma^\star_{1/a}(\cdot)$ \eqref{eq:leaky-ReLU-dual} activations, with parameter $\omega \in \Omega(d, L)$ defined in \eqref{eq:param-space-omega}.
	
	Then, the GAN estimator $\mu_{\widehat{\theta}_{m,n}}$ defined in \eqref{eq:gan-nn-B}, satisfies the following parametric rates 
	\begin{align*}
		\E d_{TV}^2\left(\nu, \mu_{\widehat{\theta}_{m,n}} \right) \precsim \sqrt{d^2 L^2 \log (dL) \left( \frac{\log m}{m} \vee \frac{\log n}{n} \right)}.
	\end{align*}
\end{theorem}

\begin{remark}
	\rm 
	The above theorem is derived using Thm.~\ref{thm:param-gans} and Cor.~\ref{cor:design}. Here we use the neural networks' architecture to perform pair regularization. Remark that our results allow for \textit{very deep} ReLU neural network with $L \precsim \sqrt{n \wedge m/\log (n \vee m) }$. The sample complexity on the dimension is $d^2$, which is desirable due to the fact that even estimating a parametric, multivariate Gaussian distribution $N(0, \Sigma)$ requires at least $d^2$ samples (see also Cor.~\ref{cor:multivariate-gaussian}).
	
	The goal here is to show with a good choice of $(\G, \F)$ suggested by the pair regularization, near-optimal sample complexity is attainable. In a nutshell, one needs to identify a pair of $(\G, \F)$ such that both $A_1(\F,\G,\nu)$ and $A_2(\G, \nu)$ are small. One sufficient choice of pair regularization to establish Thm.~\ref{thm:leaky-ReLU} and Cor.~\ref{cor:multivariate-gaussian} is: choosing the smallest $\G_\star$ with $A_2(\G, \nu) = 0$ first, and given that, selecting $\F_\star$ with $A_1(\F, \G_\star, \nu) = 0$. With the above, careful calculations of the $S_{n,m}(\F_\star,\G_\star)$ establish the results. Admittedly, we do not aim to identify the optimal pair of $(\G_\star, \F_\star)$ over the entire generator-discriminator-pair tuning domain. Such optimization can be hard. The reason is, to characterize the implicit distribution of $g_{\widehat{\theta}_{m,n}}(Z)$ given by neural networks transformations, and how it approximates general nonparametric target distribution $\nu$ is a future work outside the statistical goal of the current paper. 
\end{remark}

\begin{remark}[Relations to literature]
	
	\rm

	Investigations on the parametric rates for GANs have been considered in \cite*{bai2018approximability}, based on spectral norm-based capacity controls as regularization of networks, that is, $\forall l \in [L], \| W_l \|_{\rm op}, \| W_l^{-1} \|_{\rm op} \leq C$. The approach they took is to establish multiplicative equivalence on $d_{\F}(\mu, \nu) \asymp d_{W}(\mu, \nu)$ for $\mu, \nu \in \G$ restricted to the generator class.
	
	In contrast, we make use of the oracle inequality approach developed in an early version of the current paper \citep{liang2017well}, and the notion of pair regularization. We study through the angle of pseudo-dimensions, without requiring that the spectral radius of each $W_l, W_l^{-1}$ being bounded. This has two advantages. First, the generator class can express a wider range of densities, as we only require that $W_l$ has full rank. Second, we make explicit the dependence of the depth of the neural networks $L$ in the rate. In addition, we are able to get a better polynomial dependence on both the dimension $d$ and the depth $L$, in the error. 
\end{remark}

Finally, as a sanity check, we show that GANs can also achieve the correct dimension dependence in sample complexity, $n = O(d^2\log d)$, when estimating multivariate Gaussian with unknown mean and covariance: from classic information-theoretic lower bounds, $n = \Theta(d^2)$ samples are necessary. The example is to showcase that with the power of pair regularization, GANs can obtain provable guarantee in classic realms. 
\begin{corollary}[Multivariate Gaussian estimation]
	\label{cor:multivariate-gaussian}
	Consider $\nu \sim N(b_*, \Sigma_*)$ a multivariate Gaussian in $\mathbb{R}^d$. Consider a linear generator (neural network with no hidden layer) with input distribution $N(0, I_p)$ ($p\geq d$), and the discriminator to be a one hidden layer neural network with quadratic activation $\sigma(t) = t^2$, the GAN estimator $\mu_{\widehat{\theta}_{m,n}}$ defined in \eqref{eq:gan-nn-B}, satisfies the following rate,
	\begin{align*}
		\E d_{TV}^2\left(\nu, \mu_{\widehat{\theta}_{m,n}} \right) \precsim \sqrt{ \frac{d^2 \log d}{n} + \frac{(pd+d^2) \log (p + d)}{m} }.
	\end{align*}
\end{corollary}

\section{Discussions and Future Work}
\label{sec:discussion}

\subsection{Pair Regularization and Regularity in Optimal Transport}
We now discuss an interesting connection between the pair regularization theory developed in this paper and the regularity theory in optimal transport. To illustrate such a connection, consider the dual formulation of the Wasserstein-$2$ distance \citep[Chapter 2]{ambrosio2013user} between the input distribution $\pi$ and the target distribution $\nu$, both supported on bounded open subsets of $\mathbb{R}^d$,
\begin{align}
	\label{eq:w2-dual}
	\frac{1}{2} d_{W_2}^2 (\pi, \nu) := \sup_{f: \mathbb{R}^d \rightarrow \mathbb{R}} \int_{ \mathcal{Z}} f d\pi + \int_{\mathcal{X}} f^{c} d\nu
\end{align}
with $f^c(x):= \inf_{z \in \mathcal{Z}}  \| x- z \|^2/2 - f(z)$ is conjugate to $f$. Here the continuous function $f$ is also referred to as Kantorovich potential. Brenier showed \citep[see Theorem 2.26]{ambrosio2013user} that under mild conditions on $\pi, \nu$, the optimal $f_\star$ (up to an additive constant) that maximizes \eqref{eq:w2-dual} must satisfy:
(1) $\|z\|^2/2 - f_\star(z)$ is convex, and (2) the push-forward map $Id - \nabla f_\star:  z \mapsto x$ transforms the input distribution $\pi$ (on $\mathcal{Z}$) to the target $\nu$ (on $\mathcal{X}$)
\begin{align*}
	\nu = (Id - \nabla f_\star)_{\#} \pi \enspace.
\end{align*}
Such a transformation is exactly what is needed for learning distribution implicitly as in GANs.
In the view of GANs, the above precisely present a pair of generator and discriminator with the choice $(g_\star, f_\star)$ satisfying
\begin{align*}
	g_\star =  Id - \nabla f_\star \enspace.
\end{align*}
Conceptually, the above equation conforms with the ``pair regularization'' idea: when the generator class $\{ g_\star\}$ exhibits certain structure, the discriminator class $\{ f_\star\}$ should be picked judiciously to utilize such structure, and vice versa. This is closely related to Caffarelli's regularity theory in optimal transport \citep{caffarelli1991some, caffarelli1992regularity}: for example, when the distribution class $\pi, \nu \in \G$ (generator) satisfy that $p_\pi, p_\nu$ are $\alpha$-smooth densities in H\"{o}lder class $W^{\alpha, \infty}$, then the optimal dual potential $f\in \F$ (discriminator) is $(\alpha+2)$-smooth $W^{\alpha+2, \infty}$, under the domain conditions that $\mathcal{X}, \mathcal{Z}$ are convex with smooth boundaries.

\subsection{Statistical Advantage of GANs}

We further discuss the following question: overlooking computation, what is the advantage of GANs compared to classic nonparametric and parametric models? We use the diagram as in Fig.~\ref{fig:gan} (see Appendix) to illustrate some conclusions and conjectures.

\textit{Classic parametric models} can be viewed as the left interval (along the y-axis) in Fig.~\ref{fig:gan}, where the generator class $\G$ is limited. The discriminator is assessing how well we estimate the finite parameters, which relates to how well we are learning distributions in the parametric class. More advanced discriminator won't help. 
	
\textit{Classic nonparametric density estimation} can be viewed as the top interval (along the x-axis) in Fig.~\ref{fig:gan}. Here the evaluation metric is either $L^2$ or $L^\infty$, and by tuning the generator class $\G$ (using sieves, kernels, etc.), optimal rates are achieved when the target density lies in a certain nonparametric class. The minimax theory for the adversarial framework (Thm.~\ref{thm:optimal-rates-sobolev}) informs us, when the target is nonparametric, tuning with the generator class is optimal: there is no theoretical gain in utilizing the generator-discriminator-pair to tune. Though, with simpler evaluation metrics, one can obtain faster rates, shown in Thm.~\ref{thm:optimal-rates-sobolev}.
	
\textit{Empirical distribution}, or data memorization can be viewed as the right interval (along the y-axis) in Fig.~\ref{fig:gan}. Here the generator class is flexible enough to memorize the training data, and one should try to avoid this through regularization (Thm.~\ref{thm:nonparam-gans}).
	
For a certain target distribution $\nu$ (in between parametric and nonparametric for many realistic cases), we \textit{conjecture} that \textit{tuning with the generator-discriminator-pair} $(\G_\star, \F_\star)$ could potentially explain the empirical success of GANs on the statistical side. One can tune the generator and discriminator pair with deep neural networks, thus navigating in the two-dimensional domain balancing $A_1(\F, \G, \nu)$, $A_2(\G, \nu)$, and $S_{n, m}(\F, \G)$ simultaneously. As seen in Thm.~\ref{thm:leaky-ReLU}, the discriminator and generator classes should be chosen as a pair with matching complexities: one being too complex the other being simple does not help. In view of the optimal transport theory, the dual potential (discriminator $f_\star$), and the push-forward map (generator transformation $I_d - \nabla f_\star$) should be parametrized and chosen in a way with matching complexities. Though in a high level, the above message is important for the practical architecture design of GANs.

Admittedly, to fully understand pair regularization, one may need to rethink the class of distributions of interest. For instance, what constitutes ``low complexity'' or ``structured'' class beyond the ``smoothness'' considered in the nonparametric literature.
In this paper, we only study the statistical framework of how well GANs learn distributions, assuming the optimization, say \eqref{eq:gan-nn-B}, can be done to sufficient accuracy. 
Admittedly, computation in GANs is a considerably harder question \citep*{mescheder2017numerics, daskalakis2017training,liang2018interaction,arbel2018gradient,lucic2017gans}, which we leave as future work.

\section{Proof of Main Results}
\label{sec:proof}

\subsection{Oracle Inequalities}

We develop the oracle inequalities, which are the main tool for analyzing the rates of GANs. Remark that these are deterministic inequalities that hold generally, which could be of independent interest for further research on GAN, GMM, and SMM.

\begin{lemma}[Simple oracle inequality]
	\label{lem:oracle-ineq}
	Under the condition that $\F_D$ is a symmetric class with $\F_D = -\F_D$, the GAN estimator in \eqref{eq:gan-nonparam} satisfies
	\begin{align*}
		d_{\F_D}(\nu, \mu_n) \leq  \min_{\mu \in D_G} d_{\F_D}(\mu, \nu) + 2d_{\F_D}(\nu, \nu_n),
	\end{align*}
	where we refer the first term as the approximation error, and second as the stochastic error.
\end{lemma}

\begin{proof}[Proof of Lemma~\ref{lem:oracle-ineq}]
	For any $\mu \in \mu_G$, we know that due to the optimality of GAN in \eqref{eq:gan-nonparam},
	\begin{align*}
		d_{\F_D}(\mu, \nu_n) - d_{\F_D}(\mu_n, \nu_n) \geq 0.
	\end{align*}
	Due to the triangle inequality of IPM, we have
	\begin{align*}
		d_{\F_D}(\mu_n, \nu) & \leq d_{\F_D}(\mu_n, \nu_n) +  d_{\F_D}(\nu_n, \nu)  \\
		& \leq  d_{\F_D}(\mu, \nu_n) +  d_{\F_D}(\nu_n, \nu) \quad \text{(optimality of $\nu_n$)} \\
		& \leq  d_{\F_D}(\mu, \nu) +  d_{\F_D}(\nu, \nu_n) + d_{\F_D}(\nu_n, \nu).
	\end{align*}
	Now take $\mu = \argmin_{\mu \in \mu_G} d_{\F_D}(\mu, \nu)$, and recall that $\F_D$ is symmetric, we have
	$
		d_{\F_D}(\mu_n, \nu) \leq \min_{\mu \in \mu_G} d_{\F_D}(\mu, \nu) + 2 d_{\F_D}(\nu, \nu_n).
	$
	
\end{proof}

\begin{proof}[Proof of Lemma~\ref{lem:oracle-ineq-general}]
	For ease of notation, we abbreviate $\widehat{\theta}_{m,n}$ as $\widehat{\theta}$ in this proof when there is no confusion. Recall the GANs estimator \eqref{eq:gan-nn}, and the definition of $d_{\F}\left( \mu_{\widehat{\theta}_{m,n}}, \nu  \right)$, we have 
	\begin{align*}
		& d_{\F}\left( \mu_{\widehat{\theta}_{m,n}}, \nu  \right) = \sup_{f_\omega \in \F} \left\{ \E f_\omega \circ g_{\widehat{\theta}}(Z) - \E f_\omega (X) \right\} \\
		&\leq \sup_{f_\omega \in \F} \left\{ \E f_\omega \circ g_{\widehat{\theta}}(Z) - \widehat{\E}_n f_\omega (X) \right\} + \sup_{f_\omega \in \F} \left\{ \widehat{\E}_n f_\omega (X)- \E f_\omega (X) \right\}  \\
		&\leq \sup_{f_\omega \in \F} \left\{ \widehat{\E}_m f_\omega \circ g_{\widehat{\theta}}(Z) - \widehat{\E}_n f_\omega (X) \right\} + \sup_{f_\omega \in \F} \left\{ \E f_\omega \circ g_{\widehat{\theta}}(Z) -\widehat{\E}_m f_\omega \circ g_{\widehat{\theta}}(Z) \right\} \\
		& \hspace{3cm} + \sup_{f_\omega \in \F} \left\{ \widehat{\E}_n f_\omega (X)- \E f_\omega (X) \right\}.
	\end{align*}
	Here the first inequality we insert the quantity $\widehat{\E}_n f_\omega (X)$, and the second we insert the quantity $\widehat{\E}_m f_\omega \circ g_{\widehat{\theta}}(Z)$ to the first term.
	For any $\theta$ such that $g_\theta \in \G$, recall the optimality condition of GANs estimator
	\begin{align*}
		\sup_{f_\omega \in \F} \left\{ \widehat{\E}_m f_\omega \circ g_{\widehat{\theta}_{m,n}}(Z) - \widehat{\E}_n f_\omega (X) \right\}  \leq  \sup_{f_\omega \in \F} \left\{ \widehat{\E}_m f_\omega \circ g_{\theta}(Z) - \widehat{\E}_n f_\omega (X) \right\},
	\end{align*}
	then one can proceed with (for any fixed $\theta$ with $g_\theta \in \G$)
	\begin{align*}
		&d_{\F}\left( \mu_{\widehat{\theta}_{m,n}}, \nu  \right)  \\
		& \leq \sup_{f_\omega \in \F} \left\{ \widehat{\E}_m f_\omega \circ g_{\theta}(Z) - \widehat{\E}_n f_\omega (X) \right\} \quad \text{(optimality of $\widehat{\theta}_{m,n}$)}\\
		& \quad \quad +  \sup_{f_\omega \in \F} \left\{ \E f_\omega \circ g_{\widehat{\theta}}(Z) -\widehat{\E}_m f_\omega \circ g_{\widehat{\theta}}(Z) \right\}             + \sup_{f_\omega \in \F} \left\{ \widehat{\E}_n f_\omega (X)- \E f_\omega (X) \right\} \\
		& \leq \sup_{f_\omega \in \F} \left\{ \widehat{\E}_m f_\omega \circ g_{\theta}(Z) - \E f_\omega \circ g_{\theta}(Z) \right\} + \sup_{f_\omega \in \F} \left\{\E f_\omega \circ g_{\theta}(Z) -  \E f_\omega (X) \right\} \\
		& \quad \quad + \sup_{f_\omega \in \F} \left\{  \E f_\omega (X) - \widehat{\E}_n f_\omega (X) \right\} \quad \text{(insert $\E f_\omega \circ g_{\theta}(Z)$ and $ \E f_\omega (X)$)} \\
		& \quad \quad +  \sup_{f_\omega \in \F} \left\{ \E [f_\omega \circ g_{\widehat{\theta}}(Z)] -\widehat{\E}_m [f_\omega \circ g_{\widehat{\theta}}(Z)] \right\} + \sup_{f_\omega \in \F} \left\{ \widehat{\E}_n [f_\omega (X)- \E f_\omega (X) \right\} \\
		&\leq 2 \sup_{f_\omega \in \F} \left\{ \widehat{\E}_n f_\omega (X)- \E f_\omega (X) \right\} + \sup_{f_\omega \in \F} \left\{ \widehat{\E}_m f_\omega \circ g_{\theta}(Z) - \E f_\omega \circ g_{\theta}(Z) \right\} \\
		&\quad \quad  + \sup_{f_\omega \in \F} \left\{ \E f_\omega \circ g_{\widehat{\theta}}(Z) -\widehat{\E}_m f_\omega \circ g_{\widehat{\theta}}(Z) \right\} + \sup_{f_\omega \in \F} \left\{\E f_\omega \circ g_{\theta}(Z) -  \E f_\omega (X) \right\} \;,
	\end{align*}
	where the last step uses the fact that $f_\omega \in \F$ then $-f_\omega \in \F$. As the above holds for any $\theta$ such that $g_\theta \in \G$, we know then (by moving the last term to the LHS)
	\begin{align*}
		&d_{\F}\left( \mu_{\widehat{\theta}_{m,n}}, \nu  \right) - d_{\F}(\mu_\theta, \nu)\\
		&\leq 2 d_{\F}\left( \widehat{\nu}^n, \nu \right) + d_{\F}\left( \widehat{\mu}^m_\theta, \mu_\theta \right) +\sup_{f_\omega \in \F} \left\{ \E f_\omega \circ g_{\widehat{\theta}}(Z) -\widehat{\E}_m f_\omega \circ g_{\widehat{\theta}}(Z) \right\}  \\
		& \leq 2 d_{\F}\left( \widehat{\nu}^n, \nu \right) + d_{\F}\left( \widehat{\mu}^m_\theta, \mu_\theta \right) +\sup_{f_\omega \in \F, g_\theta \in \G} \left\{ \E f_\omega \circ g_{\theta}(Z) -\widehat{\E}_m f_\omega \circ g_{\theta}(Z) \right\} \\
		& \leq 2 d_{\F}\left( \widehat{\nu}^n, \nu \right) + d_{\F}\left( \widehat{\mu}^m_\theta, \mu_\theta \right) + d_{\F \circ \G}(\widehat{\pi}^m, \pi) \;.
	\end{align*}
	Here the second inequality is using the fact that $g_{\widehat{\theta}} \in \G$.
\end{proof}

\subsection{Minimax Optimal Rates}

We start with an equivalent definition of the Sobolev class for $W^{\alpha, q}(r)$ for $q=2$ is through the coefficients of the Fourier series. 
The following is also called the Sobolev ellipsoid.
The definition (for $q=2$) naturally extends to non-integer $\alpha \in \mathbb{R}_{\geq 0}$ through the Bessel potential. Denote $\mathbf{F}[f](\xi)$ denotes the Fourier transform of $f(x)$, and $\mathbf{F}^{-1}$ as its inverse.
\begin{definition}
	\label{def:sobolev-frac}
	For $\alpha \in \mathbb{R}_{\geq 0}$, the Sobolev class $W^{\alpha, 2}(r)$ definition extends to non-integer $\alpha$,
	\begin{align*}
		W^{\alpha}(r) := \left\{ f \in \Omega \rightarrow \mathbb{R}: \left\| \mathbf{F}^{-1}\left[ (1+|\xi|^2)^{\frac{\alpha}{2}} \mathbf{F} [f](\xi) \right] \right\|_2 \leq r \right\}.
	\end{align*}
\end{definition}

\begin{definition}[Sobolev ellipsoid]
	\label{eq:sobolev-ellipsoid}
	Let $\theta = \{\theta_\xi, \xi = (\xi_1, \ldots, \xi_d) \in \mathbb{N}^d \}$ collects the coefficients of the Fourier series, define
	\begin{align*}
		\Theta^{\alpha}(r) := \left\{ \theta \in \mathbb{N}^d \rightarrow \mathbb{R} : \sum_{\xi \in \mathbb{N}^d} (1+\sum_{i=1}^d \xi_i^2 )^{\alpha} \theta^2_{\xi} \leq r^2 \right\}.
	\end{align*}
\end{definition}
It is clear that $\Theta^\alpha(r)$ (frequency domain) is an equivalent representation of $W^{\alpha}(r)$ (spatial domain, Def.~\ref{def:sobolev-frac}) in $L^2(\mathbb{N}^d)$ for trigonometric Fourier series. For more details on Sobolev classes, we refer the readers to \cite*{nemirovski2000topics, tsybakov2009introduction, nickl2007bracketing}.


\begin{proof}[Proof of Theorem~\ref{thm:optimal-rates-sobolev}]
	
The proof consists of three main parts, the upper bound and the nonparametric minimax lower bound, and the parametric lower bound. In the proof, for simplicity, we only consider $\alpha, \beta \in \mathbb{N}_{\geq 0}$. Extensions to the $\mathbb{R}_{\geq 0}$ follows the same proof idea. 

\paragraph{\it Step 1: upper bound}
Recall that the base measure $\pi$ to be a uniform measure on $[0,1]^d$ (Lebesgue measure). 
	For the density function $p_\nu(x)$ of $\nu$ w.r.t. the Lebesgue measure, we can represent it in the Fourier trigonometric series form
	\begin{align*}
		p_\nu(x)  = \sum_{\xi \in \mathbb{N}^d} \theta_{\xi}(\nu) \psi_{\xi}(x),\quad \text{$\theta(\nu) \in \mathbb{N}^d$ denotes the coefficients of $\nu$}
	\end{align*}
	with the tensorized basis $\psi_{\xi}(x) = \prod_{i=1}^d \psi_{\xi_i}(x_i)$ .
	We construct the following estimator $p_{\widetilde{\nu}_n}$, with a cut-off parameter $M$ to be determined later,
	\begin{align*}
		p_{\widetilde{\nu}_n}(x):= \sum_{\xi \in \mathbb{N}^d}  \widetilde{\theta}_{\xi}(\nu) \psi_{\xi}(x),
	\end{align*}
	where based on i.i.d. samples $X^{(1)}, X^{(2)}, \ldots X^{(n)} \sim \nu$
	\begin{equation*}
		\widetilde{\theta}_{\xi}(\nu) :=
		\begin{cases}
		\frac{1}{n} \sum_{j=1}^n \prod_{i=1}^d \psi_{\xi_i}(X^{(j)}_i), & \text{for $\xi$ satisfies}~ \|\xi\|_{\infty} \leq M \\
		0, &\text{otherwise}
		\end{cases}.
	\end{equation*}
	Note $\widetilde{\nu}_n$ filters out all the high frequency (less smooth) components, when the multi-index $\xi$ has the largest coordinate above $M$. Similarly, expand the discriminator function $f \in \F$ in the same Fourier basis,
	$$
	f(x) = \sum_{\xi \in \mathbb{N}^d} \theta_{\xi}(f) \psi_{\xi}(x).
	$$
	Recall the Sobolev ball Def.~\ref{eq:sobolev-ellipsoid}, for any $p_\nu(x) \in W^{\alpha}(r)$, we have for the estimator $\widetilde{\nu}_n$
	\begin{align*}
		\E d_{\F}(\nu, \widetilde{\nu}_n) &= \E \sup_{f \in \F} \int_{\Omega} f(x) \left(p_\nu(x) - p_{\widetilde{\nu}_n}(x) \right) dx \\
		&= \E \sup_{f \in \F} \sum_{\xi \in \mathbb{N}^d} \theta_{\xi}(f) \left( \widetilde{\theta}_{\xi}(\nu) - \theta_{\xi}(\nu) \right) \\
		&= \E \sup_{f \in \F} \left\{ \sum_{\xi \in [M]^d} \theta_{\xi}(f) \left( \widetilde{\theta}_{\xi}(\nu) - \theta_{\xi}(\nu) \right) + \sum_{\xi \in \mathbb{N}^d \backslash [M]^d} \theta_{\xi}(f) \theta_{\xi}(\nu) \right\} \\
		&\leq \E \sup_{f \in \F}  \sum_{\xi \in [M]^d} \theta_{\xi}(f) \left( \widetilde{\theta}_{\xi}(\nu) - \theta_{\xi}(\nu) \right) + \E \sup_{f \in \F} \sum_{\xi \in \mathbb{N}^d \backslash [M]^d} \theta_{\xi}(f) \theta_{\xi}(\nu).
	\end{align*}
	
	For the truncated first term, we know
		\begin{align}
			& \E \sup_{f \in \F}  \sum_{\xi \in [M]^d} \theta_{\xi}(f) \left( \widetilde{\theta}_{\xi}(\nu) - \theta_{\xi}(\nu) \right) \nonumber\\
			& \leq \E \sup_{f \in \F} \left\{\sum_{\xi \in [M]^d} ( 1+\|\xi\|_2^2 )^{\beta} \theta^2_{\xi}(f) \right\}^{1/2} \left\{\sum_{\xi \in [M]^d} ( 1+\|\xi\|_2^2 )^{-\beta} \left( \widetilde{\theta}_{\xi}(\nu) - \theta_{\xi}(\nu) \right)^2 \right\}^{1/2} \nonumber\\
			& \leq 	 \E \left\{\sum_{\xi \in [M]^d} ( 1+\|\xi\|_2^2 )^{-\beta}  \left( \widetilde{\theta}_{\xi}(\nu) - \theta_{\xi}(\nu) \right)^2 \right\}^{1/2}  \quad\text{(as
		$\sup_{f \in \F}  \sum_{\xi \in [M]^d} ( 1+\|\xi\|_2^2 )^{\beta} \theta^2_{\xi}(f) \leq 1$)} \label{eq:improvement} \\
			& \leq 	 \left\{ \sum_{\xi \in [M]^d}  ( 1+\|\xi\|_2^2 )^{-\beta}  \E \left( \widetilde{\theta}_{\xi}(\nu) - \theta_{\xi}(\nu) \right)^2 \right\}^{1/2} \quad \text{(Jensen's inequality)} \\
			& \leq 	 \sqrt{C_{d, \beta} \frac{M^{d - 2\beta} \vee 1}{n} } \nonumber
		\end{align}
	where the last line $\E \left( \widetilde{\theta}_{\xi}(\nu) - \theta_{\xi}(\nu) \right)^2 \leq  \frac{1}{n} \E_{X \sim \nu} \psi^2_{\xi}(X) \leq \frac{1}{n}$ for trigonometric series for any multi-index $\xi$. In addition, simple calculus shows that  
	$$\sum_{\xi \in [M]^d}  ( 1+\|\xi\|_2^2 )^{-\beta} \leq C'_{d, \beta} \int_{0}^{\sqrt{d}M} \frac{r^{d-1}}{(1+r^2)^\beta} dr \leq C_{d, \beta} \left(M^{d - 2\beta} \vee 1 \right).$$

	For the second term, the following inequality holds
	\begin{align*}
		&\E \sup_{f \in \F} \sum_{\xi \in \mathbb{N}^d \backslash [M]^d} \theta_{\xi}(f) \theta_{\xi}(g) \leq \E \sup_{f \in \F} \left\{\sum_{\xi \in [M]^d} \theta^2_{\xi}(f) \right\}^{1/2} \cdot \left\{\sum_{\xi \in [M]^d} \theta^2_{\xi}(g)  \right\}^{1/2} \\
		&\leq \sup_{f \in \F} \left\{ (1+M^2)^{-\beta} \sum_{\xi \in [M]^d} (1+ \| \xi\|_2^2)^{\beta}\theta^2_{\xi}(f) \right\}^{1/2} \left\{ (1+M^2)^{-\alpha} \sum_{\xi \in [M]^d} (1+ \| \xi\|_2^2)^{\alpha} \theta^2_{\xi}(g) \right\}^{1/2} \\
		&\leq r \sqrt{\frac{1}{M^{2(\alpha+\beta)}}}.
	\end{align*}
	Combining two terms, we have for any $\nu \in \G$, with the optimal choice of $M \asymp n^{\frac{1}{2\alpha+d}}$
	\begin{align}
		\sup_{\nu\in \G} \E d_{\F}(\nu, \widetilde{\nu}_n) &\leq  \inf_{M \in \mathbb{N}} \left\{  \sqrt{C \frac{M^{d-2\beta} \vee 1}{n} } + r \sqrt{\frac{1}{M^{2(\alpha+\beta)}}} \right\} \label{eq:l_a_l_b}\\
		&\precsim n^{-\frac{\alpha+\beta}{2\alpha+d}} \vee n^{-\frac{1}{2}}.  \nonumber
	\end{align}

	Let us now establish the lower bound. Again we consider the $\Omega = [0,1]^d$ as the domain, which is the same as in the upper bound.
	\paragraph{\it Step 2: nonparametric lower bound}
	
	The main idea behind the proof is to reduce the estimation problem to a multiple hypotheses testing problem that is at least as hard. 
	In this proof, it turns out the H\"{o}lder space $W^{\alpha,\infty}$---which is a subspace of the Sobolev class $W^{\alpha}$---suffices for the minimax lower bound.

	First, we need to construct multiple distributions $\nu$'s with valid densities in $W^{\alpha, \infty}(1)$.
	Specify a kernel function $K(u) = (a_1 \exp(-\frac{1}{1 - 4u^2}) - a_2) I(|u|<1/2), u \in \mathbb{R}$ for some small fixed $a_1, a_2>0$ to ensure that $K(x) \in W^{\alpha \vee \beta, \infty}(1)$, and $\int K(u) du = 0$. This is possible since $K(u) \in C^\infty$ with uniformly bounded derivatives up to order $\alpha\vee \beta$, and therefore $a_1, a_2$ are nothing but normalization factors. Let $m$ be a parameter (that depends on the sample size $n$) to be determined later, and denote $\delta_m = 1/m$. Define the hypothesis class to be (of cardinality $2^{m^d}$)
	\begin{align*}
		\Omega_\alpha = \left\{ g_w(x) =  1 + \sum_{\xi \in [m]^d} w_\xi \delta_m^{\alpha} \varphi_\xi(x) ~|~ w \in \{0,1\}^{m^d} \right\} , \\
		\Lambda_\beta = \left\{ f_v(x) = \sum_{\xi \in [m]^d} v_\xi \delta_m^{\beta} \varphi_\xi(x) ~|~ v \in \{-1,1\}^{m^d} \right\},
	\end{align*}
	where
	$\varphi_\xi(x) =  \prod_{i=1}^d K\left(\frac{x_i - \frac{\xi_i - 1/2}{m}}{\delta_m} \right)$, with $\delta_m = 1/m$.
	
	Let us verify (1) $\Omega_\alpha \subset W^{\alpha, \infty}(r)$ for some $r$, and that each element in the hypothesis set is a valid density; (2) $\Lambda_\beta \subset W^{\beta, \infty}(1)$. 
	To start, for any multi-index $\gamma$ such that $|\gamma| \leq \alpha$, $\gamma \neq 0$,
	\begin{align*}
		\| D^{(\gamma)} g_w \|_{\infty} \leq \sup_{\xi \in [m]^d} \delta_m^{\alpha} \| D^{(\gamma)} \varphi_\xi \|_{\infty} = \delta_m^{\alpha - |\gamma|}  \| D^{(\gamma)} K(u) \|_{\infty} \leq \delta_m^{\alpha - |\gamma|} \leq 1.
	\end{align*}
    Similarly for $\forall \gamma, |\gamma| \leq \beta$, we know $\| D^{(\gamma)} f_v(x) \|_{\infty} \leq \delta_m^{\beta - |\gamma|} \leq 1.$
	We also need to bound $\| g_w \|_{\infty}$, for any $w$
	\begin{align}
		\label{eq:infty-bound}
		\| g_w \|_{\infty} \leq 1 +  \delta_m^{\alpha} \sup_{\xi \in [m]^d} \| \varphi_\xi(x) \|_\infty \leq 1 +  \delta_m^{\alpha}  \leq 1 + 1/100,
	\end{align}
	and $\inf_{x} g_w (x) \geq 1 - 1/100 >0$, when $m$ is large enough. So far we have shown $\Omega_\alpha \subset W^{\alpha, \infty}(r)$ and $\Lambda_\beta \subset W^{\beta, \infty}(1)$. 
	Last, we can check that $g_\omega$ is a proper density as we know $g_w(x) \geq 0$, and 
	$\int \varphi_\xi(x) dx = \prod_{i=1}^d \int K\left(\frac{x_i - \frac{\xi_i - 1/2}{m}}{\delta_m} \right) dx_i = 0$, $\int g_{\omega}(x) dx =  1+ \sum_{\xi \in [m]^d} w_\xi \delta_m^{\alpha} \int \varphi_\xi(x) dx = 1.$

	To select hypotheses within $\Omega_\alpha$ that are hard to distinguish with finite samples, we employ the Varshamov-Gilbert construction in conjunction with Fano's inequality (we use the version in Lemma \ref{lem:fano}). The technicality is to construct multiple hypotheses that are separated w.r.t. the adversarial loss, then show that the hypotheses are close in the statistical sense with finite samples. The Varshamov-Gilbert construction (Lemma 2.9 in \cite{tsybakov2009introduction}) claims that: for any $h \in \mathbb{N}$, there exists a subset $\{ w^{(0)}, \ldots, w^{(H)} \} \subset \{0,1\}^{h}$ with cardinality $H$, such that $w^{(0)} = (0, \ldots, 0)$,  $\rho(w^{(j)}, w^{(k)}) \geq \frac{h}{8}, ~\forall~j, k\in [H],~ j\neq k$ with $\rho(w, w')$ denoting the Hamming distance between $w$ and $w'$ on the hypercube, and that $\log H \geq \frac{h}{8} \log 2$.
	In our case set $h = m^d$.
	For the loss function, any $w, w' \in \{ w^{(0)}, \ldots, w^{(H)} \}$ with $w \neq w'$ satisfies
	\begin{align*}
		d_\F(g_w, g_{w'}) &:= \sup_{f \in W^{\beta}(1)} \int f(x) g_w(x)  dx - \int f(x) g_{w'}(x) dx \\
		& \geq \sup_{f \in W^{\beta, \infty}(1)} \int f(x) g_w(x)  dx - \int f(x) g_{w'}(x) dx \\
		& \geq \sup_{f \in \Lambda_{\beta}} \int f(x) \left(g_w(x) - g_{w'}(x) \right) dx \\
		& = \sup_{v \in \{-1, +1\}^{m^d}} \delta_m^{\alpha +\beta} \sum_{\xi\in [m]^d} v_{\xi} (w_{\xi} - w_{\xi}') \int  \varphi_\xi^2(x)  dx \\
		& =  \delta_m^{\alpha +\beta + d} \sum_{\xi \in [m]^d} \mathbbm{I}(w_\xi \neq w_{\xi}') \int \prod_{i=1}^d K^2\left(u_i \right) du  \\
		&\geq c\cdot \delta_m^{\alpha +\beta + d} \rho(w, w') \geq  c\cdot \frac{m^d}{8} \delta_m^{\alpha +\beta + d} \asymp \delta_m^{\alpha +\beta}.
	\end{align*}

	Now let's show that based $n$ i.i.d. data generated from density $g_w(x)$, it is hard to distinguish the hypotheses. For any $\omega \in \{ 0, 1\}^{h}$ with $h=m^d$, it induces a valid density $g_{\omega}$ and thus a probability distribution denoted as $\mathcal{P}_{\omega}$. Use $\mathcal{P}_{\omega}^{\otimes n}$ to denote the probability distribution of $n$-i.i.d. samples jointly.
	Note that for $|t| < 1/50$, $\log(1+t) \geq t - t^2$. Recall \eqref{eq:infty-bound} we know $\| (g_w(x) - g_{0}(x))/g_w(x)\|_{\infty} \leq \frac{1/100}{1- 1/100} \leq 1/50$, and
	\begin{align*}
		d_{KL}\left(\mathcal{P}_{w^{(j)}}^{\otimes n} || \mathcal{P}_{w^{(0)}}^{\otimes n}\right) & = n \cdot d_{KL}\left(\mathcal{P}_{w^{(j)}} || \mathcal{P}_{w^{(0)}}\right) \\
		& = n \int  - \log \left( 1 + \frac{g_{0} - g_{w^{(j)}}}{g_{w^{(j)}}} \right) g_{w^{(j)}} dx \\
		& \leq n \int \frac{(g_{0} - g_{w^{(j)}})^2}{g_{w^{(j)}}} dx \leq 1.01 n \sum_{\xi \in [m]^d} \int  \delta_m^{2\alpha} \varphi_\xi^2(x)  dx \\
		& \leq 1.01 n \sum_{\xi \in [m]^d} \int  \delta_m^{2\alpha+d} \prod_{i=1}^d K^2\left(u_i \right) du
		 \precsim n \delta_m^{2\alpha +d} m^d.
	\end{align*}
	Therefore if we choose an integer $m \asymp n^{-\frac{1}{2\alpha + d}}$ and $\delta_m = 1/m$,
	$$\frac{1}{H} \sum_{j=1}^H d_{ KL}\left(\mathcal{P}_{w^{(j)}}^{\otimes n} || \mathcal{P}_{w^{(0)}}^{\otimes n} \right) \leq  c \cdot n \delta_m^{2\alpha +d} m^d = c' \cdot m^d \leq c'' \cdot \log H.$$
	Using the Fano's inequality, the lower bound for adversarial loss is of the order $\delta_m^{\alpha+\beta} = n^{-\frac{\alpha+\beta}{2\alpha+d}}$, as
	\begin{align*}
		\inf_{\widetilde{\nu}_n} \sup_{\nu \in W^{\alpha}(r)}\E d_\F\left(\widetilde{\nu}_n, \nu \right) &\geq \inf_{\hat{g}} \sup_{g \in W^{\alpha,\infty}(r)} \E \sup_{f \in W^{\beta,\infty}(1)} \int f(x)\left( \hat{g}(x) - g(x) \right) dx   \\
		&\geq \inf_{\hat{w}} \sup_{w \in \{w^{(0)},\ldots, w^{(H)}\}} \E d_\F\left(g_{\hat{w}}, g_{w}\right)\\
		&\geq c \cdot \delta_m^{\alpha +\beta} \cdot \inf_{\hat{w}} \sup_{w \in \{w^{(0)},\ldots, w^{(H)}\}} \mathbb{P}_w\left( d_\F\left(g_{\hat{w}}, g_{w}\right) \geq c' \cdot \delta_m^{\alpha +\beta} \right) \\
		& \geq c'' \cdot \delta_m^{\alpha +\beta}  \frac{\sqrt{H}}{1+\sqrt{H}} \left( 1 - 2c' - \sqrt{\frac{2c'}{\log H}} \right) \quad \text{(Lemma~\ref{lem:fano})} \\
		&\geq c'' \cdot  n^{-\frac{\alpha+\beta}{2\alpha + d}}.
	\end{align*}

	\paragraph{\it Step 3: parametric lower bound}
	The parametric rate lower bound $n^{-1/2}$ can be obtained by the following reduction to a two-point hypothesis testing problem.
	Consider the uniform measure $p_{\nu_0}(x) = 1$ for $x \in [0, 1]^d$, and 
	\begin{align*}
		p_{\nu_1}(x) = 1 + \frac{1}{\sqrt{n}} K\left(x(1) - \frac{1}{2}\right).
	\end{align*}
	One can verify both $\nu_0, \mu_1$ are valid distributions on $[0, 1]^d$ with 
	\begin{align*}
		d_{\chi^2}(\nu_1^{\otimes n}, \nu_0^{\otimes n}) = (1 + d_{\chi^2}(\nu_1, \nu_0))^n - 1 = (1+c/n)^n - 1 \leq e^c-1
	\end{align*}
	where the last line uses the fact
	\begin{align}
		d_{\chi^2}(\nu_1, \nu_0) = \frac{1}{n} \int_{-1/2}^{1/2} K^2(u) du \leq \frac{c}{n}.
	\end{align}
	Therefore, by Pinsker's inequality
	$
		d_{TV}(\nu_1^{\otimes n}, \nu_0^{\otimes n}) \leq \sqrt{d_{\chi^2}(\nu_1^{\otimes n}, \nu_0^{\otimes n})/2} \leq \sqrt{(e^c-1)/2}.
	$
	Recall the fact that the kernel $K(u) \in C^\infty$ with uniformly bounded derivatives up to order $\alpha\vee \beta$ in the domain $[-1/2, 1/2]$.
	Therefore, we know $ p_{\nu_1}(x) - p_{\nu_0}(x) = \frac{1}{\sqrt{n}} K(x(1) - \frac{1}{2}) \in W^{\alpha \vee \beta}(r/\sqrt{n})$ with some absolute constant $r>0$.
	Now it is clear that $p_{\nu_0}, p_{\nu_1} \in W^{\alpha}(r')$ for any $\alpha>0$, with some proper constant $r'> 1+r/\sqrt{n}$ independent of $n$. 
	 Hence, by the Le Cam's method (Lemma 4 in \cite{cai2015law}), for any $\widetilde{\nu}_n$
	\begin{align*}
		\sup_{\nu \in W^{\alpha}(r)} \E d_\F\left(\widetilde{\nu}_n, \nu \right) &\geq \sup_{\nu \in \{ \nu_0, \nu_1 \}} \E d_\F\left(\widetilde{\nu}_n, \nu \right) \\
		& \geq  c \cdot d_{\F}(\nu_0, \nu_1) (1 - d_{TV}(\nu_1^{\otimes n}, \nu_0^{\otimes n})) \\
		& \geq c' \cdot d_{W^{\beta}(1)}(\nu_0, \nu_1) = c' \cdot r^{-1} \frac{1}{\sqrt{n}} \cdot \int_{-1/2}^{1/2} K^2(u) du  \geq c'' \cdot \frac{1}{\sqrt{n}}
	\end{align*}
	where the last step is by choosing the discriminator function $f(x) = r^{-1} \sqrt{n} [p_{\nu_0}(x) - p_{\nu_1}(x)]$ with $f \in \F \subset W^{\alpha \vee \beta}(1) \subseteq W^{\beta}(1)$.
\end{proof}

\subsection{Rates for Neural Networks}

\begin{proof}[Proof of Theorem~\ref{thm:param-gans}]
	The proof consists of three steps. Remark in this proof, we wrote $\int$ as $\int_{\Omega}$ when there is no confusion. 
	
	\paragraph{\it Step 1: $A_1(\F, \G, \nu)$ approximation term}
	Given the distribution of $g_{\widehat{\theta}_{m,n}}(Z)$ (we abbreviate $\widehat{\theta}_{m,n}$ as $\widehat{\theta}$ in this proof), by Pinsker's inequality (Lemma~\ref{lem:pinsker-bobkov}), 
	\begin{align*}
		d_{TV}^2\left(\nu, \mu_{\widehat{\theta}}\right) \leq \frac{1}{2} d_{KL}\left(\nu || \mu_{\widehat{\theta}} \right).
	\end{align*}
	The above implies that for any $X \sim \nu$
	\begin{align*}
		4 d_{TV}^2\left(\nu, \mu_{\widehat{\theta}} \right) &\leq d_{KL}\left( \nu || \mu_{\widehat{\theta}} \right) + d_{KL}\left( \mu_{\widehat{\theta}} || \nu \right) \\
		&= \int \log \frac{p_{\nu}(x)}{p_{\mu_{\widehat{\theta}}}(x)} \left(p_\nu(x) - p_{\mu_{\widehat{\theta}}}(x) \right) dx \quad \text{(for any $f_{\omega} \in \F$)} \\
		&= \int \left(\log \frac{p_{\nu}(x)}{p_{\mu_{\widehat{\theta}}}(x)} - f_\omega(x) \right) \left( p_\nu(x) - p_{\mu_{\widehat{\theta}}}(x) \right) dx + \int f_\omega(x) \left( p_\nu(x) - p_{\mu_{\widehat{\theta}}}(x) \right) dx \\
		&\leq \int \left(\log \frac{p_{\nu}(x)}{p_{\mu_{\widehat{\theta}}}(x)} - f_\omega(x) \right) \left( p_\nu(x) - p_{\mu_{\widehat{\theta}}}(x) \right) dx + d_{\F}\left( \mu_{\widehat{\theta}}, \nu  \right) \\
		&\leq \left\| \log \frac{p_{\nu}(x)}{p_{\mu_{\widehat{\theta}}}(x)} - f_\omega(x) \right\|_{\infty} \left\| p_\nu - p_{\mu_{\widehat{\theta}}} \right\|_1 + d_{\F}\left( \mu_{\widehat{\theta}}, \nu  \right) \\
		&\leq 2 \left\| \log \frac{p_\nu}{p_{\mu_{\widehat{\theta}}}} - f_\omega \right\|_{\infty} + d_{\F}\left( \mu_{\widehat{\theta}}, \nu  \right) \;,
	\end{align*}
	where the last line is due to the fact that $\mu_{\widehat{\theta}}, \nu$ are both proper probability distributions, so $\left\| p_\nu(x) - p_{\mu_{\widehat{\theta}}}(x) \right\|_1 \leq 2$. Take $f_{\omega}$ to be the one minimizes the first term on the RHS,
	\begin{align*}
		4 d_{TV}^2\left(\nu, \mu_{\widehat{\theta}}\right) \leq 2 \inf_{f_\omega \in \F} \left\|\log \frac{p_\nu}{p_{\mu_{\widehat{\theta}}}} - f_\omega \right\|_{\infty} + d_{\F}\left( \mu_{\widehat{\theta}}, \nu  \right).
	\end{align*}
	
	\paragraph{\it Step 2: $A_2(\G, \nu)$ approximation term and oracle inequality}
	Now, let's apply the oracle approach developed in Lemma~\ref{lem:oracle-ineq-general} to $d_{\F}\left( \mu_{\widehat{\theta}}, \nu  \right)$. For any $\theta$ such that $g_\theta \in \G$, we know
	\begin{align*}
		d_{\F}\left( \mu_{\widehat{\theta}}, \nu  \right) &\leq  d_{\F}(\mu_\theta, \nu) + 2 d_{\F}\left( \widehat{\nu}^n, \nu \right) + d_{\F}\left( \widehat{\mu}^m_\theta, \mu_\theta \right) + d_{\F \circ \G}(\widehat{\pi}^m, \pi) \\
		&\leq  B \cdot d_{TV}(\mu, \nu) + 2 d_{\F}\left( \widehat{\nu}^n, \nu \right) + d_{\F}\left( \widehat{\mu}^m_\theta, \mu_\theta \right) + d_{\F \circ \G}(\widehat{\pi}^m, \pi)\\
		& \leq  B  \cdot \sqrt{ \frac{1}{4} \left[ d_{KL}\left(\mu_\theta || \nu \right) + d_{KL}\left(\nu||\mu_\theta \right) \right] } \\
		&\quad\quad +  2 d_{\F}\left( \widehat{\nu}^n, \nu \right) + d_{\F}\left( \widehat{\mu}^m_\theta, \mu_\theta \right) + d_{\F \circ \G}(\widehat{\pi}^m, \pi) \\
		&\leq  B \cdot \sqrt{ \frac{1}{4} \left\| \log \frac{p_{\mu_\theta}}{p_{\nu}} \right\|_{\infty} \|p_{\mu_\theta} -  p_{\nu}\|_1 } +  2 d_{\F}\left( \widehat{\nu}^n, \nu \right) + d_{\F}\left( \widehat{\mu}^m_\theta, \mu_\theta \right) + d_{\F \circ \G}(\widehat{\pi}^m, \pi)
	\end{align*}
	where second line uses the fact that for any $f \in \F$, $\| f \|_\infty \leq B$.
	
	\paragraph{\it Step 3: the stochastic term $S_{m,n}(\F, \G)$ by empirical processes}
	Assemble the bounds, we have for any $\theta$
	\begin{align*}
		4 d_{TV}^2\left(\nu, \mu_{\widehat{\theta}_{m,n}} \right) & \leq 2 \inf_{\omega} \left\| \log \frac{p_{\nu}}{p_{\mu_{\widehat{\theta}}}} - f_\omega \right\|_{\infty}  + B \sqrt{ \frac{1}{2} \left\| \log \frac{p_{\mu_\theta}}{p_{\nu}} \right\|_{\infty} } \\
		& \quad \quad + 2 d_{\F}\left( \widehat{\nu}^n, \nu \right) + d_{\F}\left( \widehat{\mu}^m_\theta, \mu_\theta \right) + d_{\F \circ \G}(\widehat{\pi}^m, \pi)
	\end{align*}
	
	Therefore by choosing any $\theta_\star$ that minimizes $\left\| \log \frac{\mu_\theta}{\nu} \right\|_{\infty}$ over the generator class
	\begin{align*}
		\E d_{TV}^2\left(\nu, \mu_{\widehat{\theta}_{m,n}} \right) & \leq \frac{1}{2} \E\left\{ \inf_{\omega} \left\| \log \frac{p_{\nu}}{p_{\mu_{\widehat{\theta}}}}  - f_\omega \right\|_{\infty} \right\} + \frac{B}{4\sqrt{2}} \sqrt{ \inf_\theta  \left\| \log \frac{p_{\mu_\theta}}{p_{\nu}} \right\|_{\infty} } \\
		& \quad \quad + \E \left\{ 2 d_{\F}\left( \widehat{\nu}^n, \nu \right) + d_{\F}\left( \widehat{\mu}^m_{\theta_\star}, \mu_{\theta_\star} \right) + d_{\F \circ \G}(\widehat{\pi}^m, \pi) \right\} \\
		&\leq \frac{1}{2} \sup_{\theta} \inf_{\omega} \left\| \log \frac{p_{\nu}}{p_{\mu_{\theta}}}  - f_\omega \right\|_{\infty}   + \frac{B}{4\sqrt{2}}  \inf_\theta  \left\| \log \frac{p_{\mu_\theta}}{p_{\nu}}  \right\|_{\infty}^{1/2} \\
		& \quad \quad + \E \left\{ 2 d_{\F}\left( \widehat{\nu}^n, \nu \right) + d_{\F}\left( \widehat{\mu}^m_{\theta_\star}, \mu_{\theta_\star} \right) + d_{\F \circ \G}(\widehat{\pi}^m, \pi) \right\}.
	\end{align*}
	
	Apply the symmetrization in Lemma~\ref{lem:symmetrization}, 
	\begin{align*}
		&\E \left\{ 2 d_{\F}\left( \widehat{\nu}^n, \nu \right) + d_{\F}\left( \widehat{\mu}^m_{\theta_\star}, \mu_{\theta_\star} \right) + d_{\F \circ \G}(\widehat{\pi}^m, \pi) \right\} \\
		& \leq 4 \E \cR_n \left(\F\right) + 2\E \cR_m \left( \F \right) + 2 \E \cR_{m} \left( \F \circ \G \right) \\
		& \leq C \sqrt{\text{Pdim}(\F) \left( \frac{\log m}{m} \vee \frac{\log n}{n} \right)} + C \sqrt{  \text{Pdim}(\F \circ \G) \frac{\log m}{m} },
	\end{align*}
	where the last step uses the relationship between Rademacher complexity and pseudo-dimension, derived in Lemma~\ref{lem:rad-vc}.
\end{proof}

\begin{proof}[Proof of Theorem~\ref{thm:param-gans-hellinger}]
	Due to Le Cam's inequality (Lemma 2.3 in \cite{tsybakov2009introduction}), we know 
	\begin{align*}
		d_{TV}^2\left( \nu, \mu_{\widehat{\theta}_{m,n}} \right) & \leq d_{H}^2\left(\nu, \mu_{\widehat{\theta}_{m,n}} \right) = \int \left( \sqrt{p_\nu(x)} - \sqrt{p_{\mu_{\widehat{\theta}}}(x)} \right)^2 dx \\
		&= \int \frac{\sqrt{p_\nu(x)} - \sqrt{p_{\mu_{\widehat{\theta}}}(x)}}{\sqrt{p_\nu(x)}+ \sqrt{p_{\mu_{\widehat{\theta}}}(x)}} \left(p_\nu(x) - p_{\mu_{\widehat{\theta}}}(x) \right) dx \quad \text{for any $f_{\omega} \in \F$} \\
		& \leq 2 \left\| \frac{\sqrt{p_\nu} - \sqrt{p_{\mu_{\widehat{\theta}}}}}{\sqrt{p_\nu}+ \sqrt{p_{\mu_{\widehat{\theta}}}}} - f_{\omega} \right\|_{\infty} + d_{\F}\left( \mu_{\widehat{\theta}}, \nu  \right).
	\end{align*}
	Due to the oracle inequality Lemma~\ref{lem:oracle-ineq-general}, for any $\theta$
	\begin{align*}
		d_{\F}\left( \mu_{\widehat{\theta}}, \nu  \right) &\leq  d_{\F}(\mu_\theta, \nu) + 2 d_{\F}\left( \widehat{\nu}^n, \nu \right) + d_{\F}\left( \widehat{\mu}^m_\theta, \mu_\theta \right) + d_{\F \circ \G}(\widehat{\pi}^m, \pi).
	\end{align*}
	For the first term, we can further upper bound,
	\begin{align*}
		d_{\F}(\mu_\theta, \nu) &\leq  B \cdot d_{TV}(\mu_\theta, \nu)  \leq  B \cdot \sqrt{ d_{H} \left(\mu_\theta, \nu \right) } \\
		&\leq  2B \cdot \left\| \frac{\sqrt{p_\nu} - \sqrt{p_{\mu_{\theta}}}}{\sqrt{p_\nu} + \sqrt{p_{\mu_{\theta}}}} \right\|_{\infty}
	\end{align*}
	where the last line follows because
	\begin{align*}
		\sqrt{ d_{H} \left(\mu_\theta, \nu \right) } &= \sqrt{ \int \left( \frac{\sqrt{p_\nu(x)} - \sqrt{p_{\mu_{\widehat{\theta}}}(x)}}{\sqrt{p_\nu(x)}+ \sqrt{p_{\mu_{\widehat{\theta}}}(x)}} \right)^2 \left(\sqrt{p_{\nu}(x)} + \sqrt{p_{\mu_{\theta}}(x)} \right)^2 dx } \\
		& \leq \left\| \frac{\sqrt{p_\nu} - \sqrt{p_{\mu_{\theta}}}}{\sqrt{p_\nu} + \sqrt{p_{\mu_{\theta}}}} \right\|_{\infty} \sqrt{ \int 2 \left(p_\nu(x)+p_{\mu_\theta}(x) \right) dx}.
	\end{align*}
	The rest of the proof follows exactly the same as in Thm.~\ref{thm:param-gans}.
\end{proof}

\begin{proof}[Proof of Theorem~\ref{thm:leaky-ReLU}]
	The proof proceeds in three steps.
	
	\paragraph{\it Step 1: recursive formula for generator distributions}
	
	Consider the generator network realized by a multi-layer perceptron: 
	\begin{align*}
		h_1 &= \sigma(W_1 z + b_1) \\
		& \ldots \\
		h_l &= \sigma(W_l h_{l-1} + b_l) \\
		& \ldots \\
		x   &= W_L h_{L-1} + b_L \enspace.
	\end{align*}
	Denote the parameter space of interest
	\begin{align}
		\label{eq:param-space-theta}
		\theta \in \Theta(d, L) := \{ (W_l \in \mathbb{R}^{d \times d}, b_l \in \mathbb{R}^d, 1\leq l \leq L) ~|~ {\rm rank}(W_l) = d, \forall 1\leq l\leq L  \} .
	\end{align}
	Let us denote the density function of of the random variable $h_\ell$ as $p(h_{\ell})$.
	Consider the density evolution from layer $l-1$ to layer $l$ (basic change of variables with Jacobian $\partial h_{l}/ \partial h_{l-1}$)
	\begin{align*}
		\log p(h_l) &= \log p(h_{l-1}) - \log |\det \left( \frac{\partial h_{l}}{\partial h_{l-1}} \right)| \\
		&= \log p(h_{l-1}) - \log |\det W_l| - \sum_{i=1}^d \log |\sigma'( \sigma^{-1} (h_l(i)))|.
	\end{align*} 
	Recursively apply the above equality to track the density of $X$, we have
	\begin{align*}
		\log p_{\mu_\theta}(x) &= \log p(h_{L-1}) - \log |\det W_{L}|, \quad \text{where $h_{L-1} = W_{L}^{-1}(x - b_L)$} \\
		&= \log p(h_{L-2}) - \sum_{j=L-1}^{L} \log |\det W_{j}| - \sum_{i=1}^d \log |\sigma'(\sigma^{-1}(h_{L-1}(i)))|, \\
		& \ldots  \quad \quad \text{where $h_{L-2} = W_{L-1}^{-1}(\sigma^{-1}(h_{L-1}) - b_{L-1})$}\\
		& = \log p_\mu(z) - \sum_{j=1}^L \log |\det W_{j}| - \sum_{j=1}^{L-1} \sum_{i=1}^d \log |\sigma'(\sigma^{-1}(h_{j}(i)))|, \\
		& \quad \quad \quad \text{where $z = W_1^{-1} (\sigma^{-1}(h_1) - b_1)$}.
	\end{align*}
	Now consider $\mu(z) = 1$ to be the uniform measure on $z \in [0, 1]^d$. Consider leaky ReLU activation $\sigma(t) = \max(t, at)$ for $0< a \leq 1$, then $\sigma^{-1}(t) = \min(t, t/a)$, and $\log |\sigma'(t)| = \log(a) \cdot 1_{t\leq 0}  $.

	Let's consider the realizable case when $\log p_\nu(x) = \log p_{\mu_{\theta_*}}(x)$ for some $\theta_* \in \Theta(d, L)$. Denote $m_{l} := \sigma^{-1}(h_{L-l})$, for any $1\leq l \leq L-1$. Then it follows that
	\begin{align}
		\label{eq:m}
		m_1 &= \sigma^{-1}( W_{L}^{-1} x - W_{L}^{-1} b_L) \\
		m_{l} &= \sigma^{-1}(  W^{-1}_{L-l+1} m_{l-1} - W^{-1}_{L-l+1} b_{L-l+1} ), \quad 1\leq l \leq L-1.
	\end{align}
	Therefore, the density can be written out explicitly, 
	\begin{align}
		\label{eq:realizable-density-MLP}
		\log p_{\mu_\theta}(x) &= - \sum_{j=1}^L \log |\det W_{j}| - \sum_{j=1}^{L-1} \sum_{i=1}^d \log \sigma'(m_{L-j}(i))  \\
		&= - \sum_{j=1}^L \log |\det W_{j}| - \sum_{j=1}^{L-1} \sum_{i=1}^d \log \sigma'(m_{j}(i))  
	\end{align}
	In addition, we know that for any $\theta$ and $\theta_*$, $\mu_\theta$ and $\mu_{\theta_*}$ (namely $\nu$) are absolutely continuous to each other, as $\mu_{\theta}(x)>0$ for any $x \in [0, 1]^d$.
	
	\paragraph{\it Step 2: construction of discriminator networks}
	Now consider a discriminator network which follows
	\begin{align*}
		m_1 &= \sigma^{-1} ( V_1 x + c_1 )\\
		& \ldots \\
		m_{L-1} &= \sigma^{-1} (V_{L-1} m_{L-2} + c_{L-1}) \\
		q_{\omega}(x) &:= \sum_{j=1}^{L-1} \sum_{i=1}^d \log (1/a) 1_{m_j(i) \leq 0} +c_L \enspace.
	\end{align*}
	Here the parameter set is,
	\begin{align}
		\label{eq:param-space-omega}
		\omega \in \Omega(d, L) := \{(V_l \in \mathbb{R}^{d \times d}, c_l \in \mathbb{R}^d, c_L \in \mathbb{R}, 1\leq l \leq L-1 )~|~ {\rm rank}(V_l) = d, \forall 1\leq l\leq L-1  \}.
	\end{align} 
	Choose the discriminator function $w = (w_1, w_2)$ where $w_1, w_2 \in \Omega(d, L)$
	\begin{align*}
		f_\omega(x) = q_{\omega_1}(x) - q_{\omega_2}(x).
	\end{align*}
	Then we can verify that Cor.~\ref{cor:design} follows. Recall the upper bound in Theorem~\ref{thm:leaky-ReLU}, we can see that for the choice of generator and discriminator
	\begin{align*}
		\frac{1}{2} \sup_{\theta} \inf_{\omega} \left\| \log \frac{p_\nu}{p_{\mu_{\theta}}} - f_\omega \right\|_{\infty} = 0 \\
		\frac{B}{4\sqrt{2}}  \inf_\theta  \left\| \log \frac{p_{\mu_\theta}}{p_\nu} \right\|_{\infty}^{1/2} = 0
	\end{align*}
	as $\log p_v(x)$ can be realized by $\log p_{\mu_{\theta_*}}(x)$, and that for any $\theta \in \Theta(d, L)$, there exist an $\omega \in \Omega(d, L)$ such that
	\begin{align*}
		f_\omega(x) = \log p_{\nu}(x) - \log p_{\mu_{\theta}}(x).
	\end{align*}
	
	\paragraph{\it Step 3: complexity bound}
	Recall the result in \cite*{Bartlett-etal2017_COLT} on the Vapnik-Chervonenkis dimension of feed-forward neural networks (See Lemma~\ref{lem:VC-dim-bound} with degree at most $1$ and number of pieces $p+1=2$), we know for leaky-ReLU neural networks $\F$ and $\F \circ \G$ respectively by careful counting based on the constructions in the Steps 1 and 2.
	\begin{align*}
		\text{for network $\F$}: & \quad \text{number of weights}~W_{\F} \leq 2(d^2 L + 2dL)+2, \\
		 & \quad ~\text{number of units}~U_{\F} \leq 4 dL, \\
		 & \quad ~\text{depth}~L_{\F} \leq L+2 ~; \\
		\text{for network $\F \circ \G$}:& \quad \text{number of weights}~W_{\F\circ \G} \leq W_{\F} + d^2 L\\
		& \quad ~\text{number of units}~U_{\F\circ \G} \leq U_{\F} + dL, \\
		& \quad ~\text{depth}~L_{\F\circ \G} \leq L_{\F} + L ~.
	\end{align*}
	Therefore, we have the following upper bound on VC-dimension,
	\begin{align*}
		{\rm Pdim}(\F) \asymp \text{VCdim}(\F) \leq C \cdot L_\F W_\F \log U_\F  = C d^2 L^2 \log(dL),\\
		{\rm Pdim}(\F\circ \G) \asymp \text{VCdim}(\F \circ \G)  \leq C \cdot L_{\F\circ \G} W_{\F\circ \G} \log U_{\F\circ \G} \leq  C' d^2 L^2 \log(dL). 
	\end{align*}
	Finally, by Cor.~\ref{cor:design}, we have the result proved.
\end{proof}


\acks{The author acknowledges the generous support from the NSF Career award (DMS-2042473), and the William S. Fishman Faculty Research Fund at the University of Chicago Booth School of Business. The author wishes to thank Maxim Raginsky, Chris Hansen and anonymous referees for valuable feedback. This paper was previously posted as ``How well can generative adversarial networks learn densities: A nonparametric view'' available on arXiv:1712.08244, 2017. The previous version is no longer intended for publication.}

\bibliographystyle{abbrevnat}
\bibliography{bibfile}

\begin{thebibliography}{52}
\providecommand{\natexlab}[1]{#1}
\providecommand{\url}[1]{\texttt{#1}}
\expandafter\ifx\csname urlstyle\endcsname\relax
  \providecommand{\doi}[1]{doi: #1}\else
  \providecommand{\doi}{doi: \begingroup \urlstyle{rm}\Url}\fi

\bibitem[Ambrosio and Gigli(2013)]{ambrosio2013user}
Luigi Ambrosio and Nicola Gigli.
\newblock A user's guide to optimal transport.
\newblock In \emph{Modelling and Optimisation of Flows on Networks}, pages
  1--155. {Springer}, 2013.

\bibitem[Anthony and Bartlett(2009)]{anthony2009neural}
Martin Anthony and Peter~L Bartlett.
\newblock \emph{Neural network learning: Theoretical foundations}.
\newblock cambridge university press, 2009.

\bibitem[Arbel et~al.(2018)Arbel, Sutherland, Bi{\'n}kowski, and
  Gretton]{arbel2018gradient}
Michael Arbel, Dougal~J Sutherland, Miko{\l}aj Bi{\'n}kowski, and Arthur
  Gretton.
\newblock On gradient regularizers for mmd gans.
\newblock \emph{arXiv preprint arXiv:1805.11565}, 2018.

\bibitem[Arjovsky and Bottou(2017)]{arjovsky2017towards}
Martin Arjovsky and L{\'e}on Bottou.
\newblock Towards principled methods for training generative adversarial
  networks.
\newblock \emph{arXiv preprint arXiv:1701.04862}, 2017.

\bibitem[Arjovsky et~al.(2017)Arjovsky, Chintala, and
  Bottou]{arjovsky2017wasserstein}
Martin Arjovsky, Soumith Chintala, and L{\'e}on Bottou.
\newblock Wasserstein gan.
\newblock \emph{arXiv preprint arXiv:1701.07875}, 2017.

\bibitem[Arora and Zhang(2017)]{arora2017gans}
Sanjeev Arora and Yi~Zhang.
\newblock Do gans actually learn the distribution? an empirical study.
\newblock \emph{arXiv preprint arXiv:1706.08224}, 2017.

\bibitem[Arora et~al.(2017)Arora, Ge, Liang, Ma, and
  Zhang]{arora2017generalization}
Sanjeev Arora, Rong Ge, Yingyu Liang, Tengyu Ma, and Yi~Zhang.
\newblock Generalization and equilibrium in generative adversarial nets (gans).
\newblock \emph{arXiv preprint arXiv:1703.00573}, 2017.

\bibitem[Athey et~al.(2019)Athey, Imbens, Metzger, and Munro]{athey2019using}
Susan Athey, Guido~W Imbens, Jonas Metzger, and Evan~M Munro.
\newblock Using wasserstein generative adversarial networks for the design of
  monte carlo simulations.
\newblock Technical report, National Bureau of Economic Research, 2019.

\bibitem[Back and Brown(1993)]{back1993ImpliedProbabilities}
Kerry Back and David~P. Brown.
\newblock Implied {{Probabilities}} in {{GMM Estimators}}.
\newblock \emph{Econometrica}, 61\penalty0 (4):\penalty0 971--975, 1993.
\newblock ISSN 0012-9682.
\newblock \doi{10.2307/2951771}.

\bibitem[Bai et~al.(2018)Bai, Ma, and Risteski]{bai2018approximability}
Yu~Bai, Tengyu Ma, and Andrej Risteski.
\newblock Approximability of discriminators implies diversity in gans.
\newblock \emph{arXiv preprint arXiv:1806.10586}, 2018.

\bibitem[Bartlett et~al.(2017)Bartlett, Harvey, Liaw, and
  Mehrabian]{Bartlett-etal2017_COLT}
Peter~L. Bartlett, Nick Harvey, Christopher Liaw, and Abbas Mehrabian.
\newblock Nearly-tight vc-dimension bounds for piecewise linear neural
  networks.
\newblock In \emph{Proceedings of the 22nd Annual Conference on Learning Theory
  (COLT 2017)}, 2017.

\bibitem[Brown(1986)]{brown1986fundamentals}
Lawrence~D Brown.
\newblock Fundamentals of statistical exponential families: with applications
  in statistical decision theory.
\newblock Ims, 1986.

\bibitem[Caffarelli(1991)]{caffarelli1991some}
Luis~A Caffarelli.
\newblock Some regularity properties of solutions of {{Monge Ampere}} equation.
\newblock \emph{Communications on pure and applied mathematics}, 44\penalty0
  (8-9):\penalty0 965--969, 1991.

\bibitem[Caffarelli(1992)]{caffarelli1992regularity}
Luis~A Caffarelli.
\newblock The regularity of mappings with a convex potential.
\newblock \emph{Journal of the American Mathematical Society}, 5\penalty0
  (1):\penalty0 99--104, 1992.

\bibitem[Cai et~al.(2015)Cai, Liang, and Zhou]{cai2015law}
T.~Tony Cai, Tengyuan Liang, and Harrison~H. Zhou.
\newblock Law of log determinant of sample covariance matrix and optimal
  estimation of differential entropy for high-dimensional gaussian
  distributions.
\newblock \emph{Journal of Multivariate Analysis}, 137:\penalty0 161 -- 172,
  2015.

\bibitem[Canas and Rosasco(2012)]{canas2012learning}
Guillermo Canas and Lorenzo Rosasco.
\newblock Learning probability measures with respect to optimal transport
  metrics.
\newblock In \emph{Advances in Neural Information Processing Systems}, pages
  2492--2500, 2012.

\bibitem[Caponnetto and De~Vito(2007)]{caponnetto2007optimal}
Andrea Caponnetto and Ernesto De~Vito.
\newblock Optimal rates for the regularized least-squares algorithm.
\newblock \emph{Foundations of Computational Mathematics}, 7\penalty0
  (3):\penalty0 331--368, 2007.

\bibitem[Chen et~al.(2020)Chen, Liao, Zha, and Zhao]{chen2020statistical}
Minshuo Chen, Wenjing Liao, Hongyuan Zha, and Tuo Zhao.
\newblock Statistical guarantees of generative adversarial networks for
  distribution estimation.
\newblock \emph{arXiv preprint arXiv:2002.03938}, 2020.

\bibitem[Daskalakis et~al.(2017)Daskalakis, Ilyas, Syrgkanis, and
  Zeng]{daskalakis2017training}
Constantinos Daskalakis, Andrew Ilyas, Vasilis Syrgkanis, and Haoyang Zeng.
\newblock Training gans with optimism.
\newblock \emph{arXiv preprint arXiv:1711.00141}, 2017.

\bibitem[De~Vito et~al.(2019)De~Vito, M{\"u}cke, and
  Rosasco]{devito2019ReproducingKernel}
Ernesto De~Vito, Nicole M{\"u}cke, and Lorenzo Rosasco.
\newblock Reproducing kernel {{Hilbert}} spaces on manifolds: {{Sobolev}} and
  {{Diffusion}} spaces.
\newblock \emph{arXiv:1905.10913 [cs, math, stat]}, May 2019.

\bibitem[Dziugaite et~al.(2015)Dziugaite, Roy, and
  Ghahramani]{dziugaite2015training}
Gintare~Karolina Dziugaite, Daniel~M Roy, and Zoubin Ghahramani.
\newblock Training generative neural networks via maximum mean discrepancy
  optimization.
\newblock \emph{arXiv preprint arXiv:1505.03906}, 2015.

\bibitem[Farrell et~al.(2021)Farrell, Liang, and Misra]{farrell2018DeepNeural}
Max~H Farrell, Tengyuan Liang, and Sanjog Misra.
\newblock Deep neural networks for estimation and inference.
\newblock \emph{Econometrica}, 89\penalty0 (1):\penalty0 181--213, January
  2021.
\newblock \doi{10.3982/ECTA16901}.

\bibitem[Goodfellow et~al.(2014)Goodfellow, Pouget-Abadie, Mirza, Xu,
  Warde-Farley, Ozair, Courville, and Bengio]{goodfellow2014generative}
Ian Goodfellow, Jean Pouget-Abadie, Mehdi Mirza, Bing Xu, David Warde-Farley,
  Sherjil Ozair, Aaron Courville, and Yoshua Bengio.
\newblock Generative adversarial nets.
\newblock In \emph{Advances in neural information processing systems}, pages
  2672--2680, 2014.

\bibitem[Hansen(1982)]{hansen1982LargeSample}
Lars~Peter Hansen.
\newblock Large {{Sample Properties}} of {{Generalized Method}} of {{Moments
  Estimators}}.
\newblock \emph{Econometrica}, 50\penalty0 (4):\penalty0 1029--1054, 1982.
\newblock ISSN 0012-9682.
\newblock \doi{10.2307/1912775}.

\bibitem[Imbens et~al.(1995)Imbens, Johnson, and Spady]{imbens1995information}
Guido~W Imbens, Phillip Johnson, and Richard~H Spady.
\newblock Information theoretic approaches to inference in moment condition
  models.
\newblock Technical report, National Bureau of Economic Research, 1995.

\bibitem[Lei et~al.(2019)Lei, Lee, Dimakis, and Daskalakis]{lei2019sgd}
Qi~Lei, Jason~D Lee, Alexandros~G Dimakis, and Constantinos Daskalakis.
\newblock Sgd learns one-layer networks in wgans.
\newblock \emph{arXiv preprint arXiv:1910.07030}, 2019.

\bibitem[Li et~al.(2015)Li, Swersky, and Zemel]{li2015generative}
Yujia Li, Kevin Swersky, and Rich Zemel.
\newblock Generative moment matching networks.
\newblock In \emph{Proceedings of the 32nd International Conference on Machine
  Learning (ICML-15)}, pages 1718--1727, 2015.

\bibitem[Liang(2017)]{liang2017well}
Tengyuan Liang.
\newblock How well can generative adversarial networks (gan) learn densities: A
  nonparametric view.
\newblock \emph{arXiv preprint arXiv:1712.08244}, 2017.

\bibitem[Liang(2019)]{liang2019EstimatingCertain}
Tengyuan Liang.
\newblock Estimating certain integral probability metric ({{IPM}}) is as hard
  as estimating under the {{IPM}}.
\newblock \emph{arXiv preprint arXiv:1911.00730}, November 2019.

\bibitem[Liang and Stokes(2019)]{liang2018interaction}
Tengyuan Liang and James Stokes.
\newblock Interaction matters: {{A}} note on non-asymptotic local convergence
  of generative adversarial networks.
\newblock In Kamalika Chaudhuri and Masashi Sugiyama, editors, \emph{The 22nd
  International Conference on Artificial Intelligence and Statistics},
  volume~89 of \emph{Proceedings of Machine Learning Research}, pages 907--915.
  {PMLR}, April 2019.

\bibitem[Liu and Chaudhuri(2018)]{liu2018inductive}
Shuang Liu and Kamalika Chaudhuri.
\newblock The inductive bias of restricted f-gans.
\newblock \emph{arXiv preprint arXiv:1809.04542}, 2018.

\bibitem[Liu et~al.(2017)Liu, Bousquet, and Chaudhuri]{liu2017approximation}
Shuang Liu, Olivier Bousquet, and Kamalika Chaudhuri.
\newblock Approximation and convergence properties of generative adversarial
  learning.
\newblock \emph{arXiv preprint arXiv:1705.08991}, 2017.

\bibitem[Lucic et~al.(2017)Lucic, Kurach, Michalski, Gelly, and
  Bousquet]{lucic2017gans}
Mario Lucic, Karol Kurach, Marcin Michalski, Sylvain Gelly, and Olivier
  Bousquet.
\newblock Are gans created equal? a large-scale study.
\newblock \emph{arXiv preprint arXiv:1711.10337}, 2017.

\bibitem[Mair and Ruymgaart(1996)]{mair1996statistical}
Bernard~A Mair and Frits~H Ruymgaart.
\newblock Statistical inverse estimation in hilbert scales.
\newblock \emph{SIAM Journal on Applied Mathematics}, 56\penalty0 (5):\penalty0
  1424--1444, 1996.

\bibitem[McFadden(1989)]{mcfadden1989MethodSimulated}
Daniel McFadden.
\newblock A {{Method}} of {{Simulated Moments}} for {{Estimation}} of
  {{Discrete Response Models Without Numerical Integration}}.
\newblock \emph{Econometrica}, 57\penalty0 (5):\penalty0 995--1026, 1989.
\newblock ISSN 0012-9682.
\newblock \doi{10.2307/1913621}.

\bibitem[Mescheder et~al.(2017)Mescheder, Nowozin, and
  Geiger]{mescheder2017numerics}
Lars Mescheder, Sebastian Nowozin, and Andreas Geiger.
\newblock The numerics of gans.
\newblock \emph{arXiv preprint arXiv:1705.10461}, 2017.

\bibitem[Mescheder et~al.(2018)Mescheder, Geiger, and
  Nowozin]{mescheder2018training}
Lars Mescheder, Andreas Geiger, and Sebastian Nowozin.
\newblock Which training methods for gans do actually converge?
\newblock In \emph{International Conference on Machine Learning}, pages
  3478--3487, 2018.

\bibitem[Mroueh et~al.(2017)Mroueh, Li, Sercu, Raj, and
  Cheng]{mroueh2017sobolev}
Youssef Mroueh, Chun-Liang Li, Tom Sercu, Anant Raj, and Yu~Cheng.
\newblock Sobolev gan.
\newblock \emph{arXiv preprint arXiv:1711.04894}, 2017.

\bibitem[Nemirovski(2000)]{nemirovski2000topics}
Arkadi Nemirovski.
\newblock Topics in non-parametric.
\newblock \emph{Ecole d'Et{\'e} de Probabilit{\'e}s de Saint-Flour},
  28:\penalty0 85, 2000.

\bibitem[Nickl and P{\"o}tscher(2007)]{nickl2007bracketing}
Richard Nickl and Benedikt~M P{\"o}tscher.
\newblock Bracketing metric entropy rates and empirical central limit theorems
  for function classes of besov-and sobolev-type.
\newblock \emph{Journal of Theoretical Probability}, 20\penalty0 (2):\penalty0
  177--199, 2007.

\bibitem[Pakes and Pollard(1989)]{pakes1989SimulationAsymptotics}
Ariel Pakes and David Pollard.
\newblock Simulation and the {{Asymptotics}} of {{Optimization Estimators}}.
\newblock \emph{Econometrica}, 57\penalty0 (5):\penalty0 1027--1057, 1989.
\newblock ISSN 0012-9682.
\newblock \doi{10.2307/1913622}.

\bibitem[Pollard(1990)]{pollard1990empirical}
David Pollard.
\newblock Empirical processes: theory and applications.
\newblock 1990.

\bibitem[Singh and P{\'o}czos(2018)]{singh2018minimax}
Shashank Singh and Barnab{\'a}s P{\'o}czos.
\newblock Minimax distribution estimation in wasserstein distance.
\newblock \emph{arXiv preprint arXiv:1802.08855}, 2018.

\bibitem[Singh et~al.(2018)Singh, Uppal, Li, Li, Zaheer, and
  P{\'o}czos]{singh2018nonparametric}
Shashank Singh, Ananya Uppal, Boyue Li, Chun-Liang Li, Manzil Zaheer, and
  Barnab{\'a}s P{\'o}czos.
\newblock Nonparametric density estimation under adversarial losses.
\newblock \emph{arXiv preprint arXiv:1805.08836}, 2018.

\bibitem[S{\o}nderby et~al.(2016)S{\o}nderby, Caballero, Theis, Shi, and
  Husz{\'a}r]{sonderby2016amortised}
Casper~Kaae S{\o}nderby, Jose Caballero, Lucas Theis, Wenzhe Shi, and Ferenc
  Husz{\'a}r.
\newblock Amortised map inference for image super-resolution.
\newblock \emph{arXiv preprint arXiv:1610.04490}, 2016.

\bibitem[Stone(1982)]{stone1982optimal}
Charles~J Stone.
\newblock Optimal global rates of convergence for nonparametric regression.
\newblock \emph{The annals of statistics}, pages 1040--1053, 1982.

\bibitem[Tsybakov(2009)]{tsybakov2009introduction}
Alexandre~B Tsybakov.
\newblock \emph{Introduction to nonparametric estimation}.
\newblock Springer Series in Statistics. Springer, New York, 2009.

\bibitem[van Handel(2014)]{van2014probability}
Ramon van Handel.
\newblock Probability in high dimension.
\newblock Technical report, PRINCETON UNIV NJ, 2014.

\bibitem[Wassermann(2006)]{wassermann2006all}
Larry Wassermann.
\newblock \emph{All of nonparametric statistics}.
\newblock Springer Science+ Business Media, New York, 2006.

\bibitem[Weed and Berthet(2019)]{weed2019estimation}
Jonathan Weed and Quentin Berthet.
\newblock Estimation of smooth densities in wasserstein distance.
\newblock \emph{arXiv preprint arXiv:1902.01778}, 2019.

\bibitem[Yarotsky(2017)]{yarotsky2016error}
Dmitry Yarotsky.
\newblock Error bounds for approximations with deep relu networks.
\newblock \emph{Neural Networks}, 94:\penalty0 103--114, 2017.

\bibitem[Yu(2013)]{yu2013stability}
Bin Yu.
\newblock Stability.
\newblock \emph{Bernoulli}, 19\penalty0 (4):\penalty0 1484--1500, 2013.

\end{thebibliography}


\newpage

\appendix

\begin{figure}[ht!]
\centering
\includegraphics[width=0.6\textwidth]{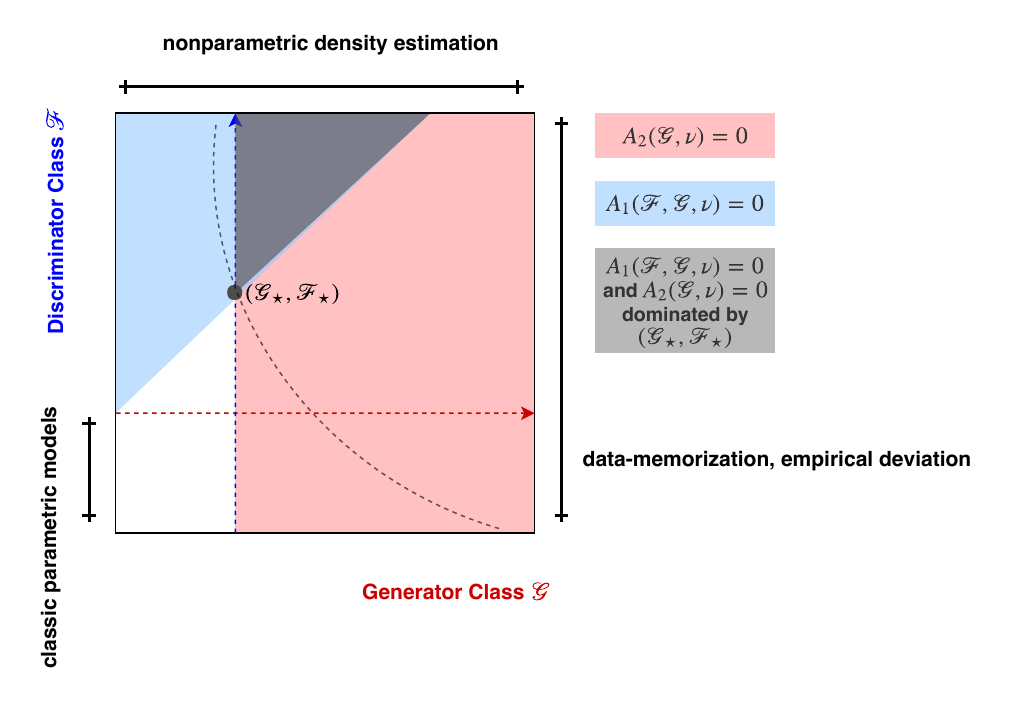}
\caption{Diagram for generator-discriminator-pair regularization.}
\label{fig:gan}
\end{figure}

\section{Remaining Proofs}
\label{sec:append}
\subsection{Other Theorems and Corollaries}

\begin{proof}[Proof of Theorem~\ref{thm:optimal-rates-rkhs}]
	The proof logic of this corollary follows similarly as in Theorem~\ref{thm:optimal-rates-sobolev}. We need to adapt the proof to the density ratio w.r.t. the general base measure $\pi$. Express $f \in \F$ under the eigenfunctions
	\begin{align*}
		f(x) = \sum_{i \in \mathbb{N}} f_i \psi_i(x), ~\text{with}~ \sum_i t_i^{-1} f_i^2 \leq 1
	\end{align*}
	where $t_i \asymp i^{-\kappa}$ and $f_i = \int f\psi_i  d\pi$ are the coefficients. Consider the series representation of the target density $d \nu /d\pi$ w.r.t. the base measure $\pi$
	\begin{align*}
		\frac{d\nu}{d\pi}(x) =  \sum_{i \in \mathbb{N}} \nu_i \psi_i(x), ~\text{then}\\
		\| \cT_{\pi}^{-(\lambda-1)/2} \frac{d \nu}{d\pi} \|_{\cH} \leq r ~\text{is equivalent to}~ \sum_i t_i^{-\lambda} \nu_i^2 \leq r^2 .
	\end{align*}
	Define the regularized density
	\begin{align*}
		\frac{d \widetilde{\nu}_n}{d\pi}(x):= \sum_{i \in \mathbb{N}}  \widetilde{\nu}_i \psi_{i}(x),
	\end{align*}
	where based on i.i.d. samples $X^{(1)}, X^{(2)}, \ldots X^{(n)} \sim \nu$
	\begin{equation*}
		\widetilde{\nu}_i :=
		\begin{cases}
		\frac{1}{n} \sum_{j=1}^n \psi_{i}(X^{(j)}), & \text{for}~ i \leq M \\
		0, &\text{otherwise}
		\end{cases}.
	\end{equation*}
	Follow the sample logic as in the proof of Theorem~\ref{thm:optimal-rates-sobolev}, we have for any $\nu(x) \in \G$, with the optimal choice of an integer $M \asymp n^{\frac{1}{\lambda \kappa+1}}$, the following holds
	\begin{align*}
		\E d_{\F}(\nu, \widetilde{\nu}_n) & = \E \int f (d\nu - d\widetilde{\nu}_n) \\
		&= \E \int f \left( \frac{d\nu}{d\pi} - \frac{d \widetilde{\nu}_n}{d\pi} \right) d\pi \\
		&\leq \E \sup_{f \in \F}  \sum_{i\leq M} f_i \left(   \nu_i - \widetilde{\nu}_i \right) +  \E \sup_{f \in \F} \sum_{i > M} f_i  \nu_i.\\
		& \leq  \sqrt{\sum_{i\leq M} t_i^{-1} f_i^2 } \sqrt{\sum_{i\leq M} t_i \E (\widetilde{\nu}_i - \nu_i)^2 } +  Cr t_M^{\frac{\lambda+1}{2}} \\
		& \asymp \inf_{M \in \mathbb{N}} \left\{  \sqrt{C \frac{M^{1-\kappa} \vee 1}{n} } +  C r \sqrt{\frac{1}{M^{\kappa(\lambda+1)}}} \right\} \\
		&\precsim n^{-\frac{(\lambda+1)\kappa}{2\lambda \kappa+2}} \vee n^{-\frac{1}{2}}.  \nonumber
	\end{align*}
	
\end{proof}

\begin{proof}[Proof of Theorem~\ref{thm:nonparam-gans}]
	Consider first the Wasserstein GAN case.
	By the entropy integral Lemma~\ref{lem:symmetrization}, 
	 if $\F_D$ consists of $L$-Lipschitz functions (Wasserstein GAN) on $\mathbb{R}^d$, $d\geq 2$, plug in the $\ell_\infty$-covering number bound for Lipschitz functions,
		\begin{align*}
			&\log \mathcal{N}(\epsilon, \F_{D}, \| \cdot \|_{\infty}) \leq C \left( \frac{L}{\epsilon} \right)^d, \\
			&\E d_{\F_D}(\nu, \widehat{\nu}^{n}) \leq 2 \inf_{0<\delta<1/2} \left( 4\delta + \frac{8 \sqrt{2}}{\sqrt{n}} \int_{\delta}^{1/2}  \sqrt{\log \mathcal{N}(\epsilon, \F_{D}, \| \cdot \|_{\infty}) } d \epsilon \right) \\
			&  \quad \quad \leq 16 \left( \frac{4\sqrt{2C}}{d-2}\right)^{\frac{2}{d}} L n^{-\frac{1}{d}} = \mathcal{O}\left( \left(\frac{C}{d^2 n} \right)^{-\frac{1}{d}} \right).
		\end{align*}
		This matches the best known bound as in \cite{canas2012learning} (Section 2.1.1).

	Let's consider now the Sobolev GAN when $\F_D$ denotes Sobolev class $W^{\beta, 2}$ on $\mathbb{R}^d$. Recall the entropy number estimate for $W^{\beta, 2}$ \citep{nickl2007bracketing}, we have
		\begin{align*}
			&\log \mathcal{N}(\epsilon, \F_{D}, \| \cdot \|_{\infty}) \leq C \left( \frac{1}{\epsilon} \right)^{\frac{d}{\beta} \vee 2},\\
			&\E d_{\F_D}(\nu, \widehat{\nu}^{n}) \leq \mathcal{O}\left( n^{-\frac{\beta}{d}} + \frac{\log n}{\sqrt{n}} \right).
		\end{align*}
		
		For the regularized distribution as plug-in, one can apply Lemma~\ref{lem:oracle-ineq} and Theorem~\ref{thm:optimal-rates-sobolev} to obtain the claimed result.
\end{proof}

\begin{proof}[Proof of Corollary~\ref{coro:GMM}]
	First, let us show why GMM is a special case of the adversarial framework. Denote a vectored-value function $\Phi(\cdot): \Omega  \rightarrow \mathbb{R}^{K}$ with elements $\Phi(x)[k] = \phi_k(x)$. \eqref{eq:gan-nonparam} is equivalent to 
	\begin{align}
		\min_{\mu \in \D_G} \big(\E_{Y\sim \mu} \Phi(Y) - \E_{X\sim \widehat\nu^n} \Phi(Y) \big)^\top \bW \big(\E_{Y\sim \mu} \Phi(Y) - \E_{X\sim \widehat\nu^n} \Phi(Y) \big)
	\end{align}
	which is precisely moment-matching between $\mu$ and $\widehat\nu^n$ with weight matrix $\bW$. By Lemma~\ref{lem:oracle-ineq}, we need to bound $\E d_{\F_D}(\nu, \widehat\nu^n)$ now. Using the symmetrization Lemma~\ref{lem:symmetrization}, $\E d_{\F_D}(\nu, \widehat\nu^n) \leq 2\E \cR_n \left(\F_D\right)$. Let's calculate the Rademacher complexity 
	\begin{align*}
		\E \cR_n \left(\F_D\right) &= \frac{1}{n} \E_{X_1,\ldots, X_n\sim \nu} \E_{\epsilon} \sup_{\omega: \omega^\top \bW^{-1} \omega \leq 1}\sum_{l \in [n]} \epsilon_l  \sum_{i \in [K]} \omega_i \phi_{i}(X_l) \\
		&= \frac{1}{n} \E_{X_1,\ldots, X_n\sim \nu} \sqrt{ \sum_{i,j\in [K]} [\sum_{l\in [n]} \epsilon_l \phi_i(X_l)]\bW_{ij}[\sum_{l\in [n]} \epsilon_l \phi_j(X_l)] } \\
		& \leq \frac{1}{n} \E_{X_1,\ldots, X_n\sim \nu}\sqrt{ \E_{\epsilon}  \sum_{i,j\in [K]} [\sum_{l\in [n]} \epsilon_l \phi_i(X_l)]\bW_{ij}[\sum_{l\in [n]} \epsilon_l \phi_j(X_l)] } \\
		& = \frac{1}{n} \E_{X_1,\ldots, X_n\sim \nu}\sqrt{ \sum_{i,j\in [K]} \sum_{l\in [n]} \phi_i(X_l)\bW_{ij} \phi_j(X_l) } \\
		& \leq \sqrt{\frac{\E_{X\sim \nu}[\sum_{i,j\in [K]} \bW_{ij} \phi_i(X) \phi_j(X) ]}{n}}
	\end{align*}
	where the third and the fourth steps uses the Jensen's inequality.
\end{proof}

\begin{proof}[Proof of Corollary~\ref{cor:wass}]
	
	Now let's consider Wasserstein distance. Consider in addition the Lipschitz constants of $\F$ to be $L_{\F}$, and $\G$ to be $L_{\G}$, namely
	\begin{align*}
		| f_{\omega}(x) - f_{\omega}(x') | \leq L_{\F} \| x - x' \| \\
		\| g_{\theta} (z) -  g_{\theta} (z') \| \leq L_{\G} \| z - z' \|
	\end{align*}
	
	Consider first the case when $Z\sim N(0, I_d)$ (unbounded). Then for any $f \in Lip(1)$, we know
	\begin{align}
		f(g_\theta(z)) \in Lip(L_{\G}).
	\end{align}
	In other words, $f \circ g_{\theta} (Z)$ are $L_{\G}^2$ sub-Gaussian (Lemma~\ref{lem:pinsker-bobkov}), 
	therefore
	\begin{align*}
		d_{W}^2 \left( \nu, \mu_{\widehat{\theta}} \right) \leq 2 L_{\G}^2 \cdot  d_{KL} \left( \nu || \mu_{\widehat{\theta}} \right) 
	\end{align*}
	and
	\begin{align*}
		d_{\F} \left( \nu, \mu_{\theta} \right) &\leq L_{\F}  \cdot d_{W}\left( \nu, \mu_{\theta} \right) \\
		& \leq \sqrt{2} L_{\F} L_{\G} \sqrt{ d_{KL} \left( \nu || \mu_{\theta} \right)  }.
	\end{align*}
	Follow the analysis with as in the TV distance, we have
	\begin{align*}
		\E d_W^2 \left( \nu, \mu_{\widehat{\theta}} \right) & \leq  L_{\G}^2 \sup_{\theta} \inf_{\omega} \left\| \log \frac{p_\nu}{p_{\mu_{\theta}}} - f_\omega \right\|_{\infty} +   L_{\G}^3 L_{\F} \inf_\theta  \left\| \log \frac{p_{\mu_\theta}}{p_\nu} \right\|_{\infty}^{1/2} \\
		& \quad \quad + C \sqrt{\text{Pdim}(\F) \left( \frac{\log m}{m} \vee \frac{\log n}{n} \right)} + C \sqrt{  \text{Pdim}(\F \circ \G) \frac{\log m}{m} }.
	\end{align*}
	
	Consider then the case when $z,x \in [0, 1]^d$ in bounded region, we know
	\begin{align}
		\| g_{\theta} (z) -  g_{\theta} (z') \| \leq L_{\G} \sqrt{d}
	\end{align}
	Therefore $\| g_{\theta}(z) \| \leq M + L_{\G} \sqrt{d}$, and the support of $g_{\theta}(Z)$ lies in an $\ell_2$ ball with $R:= M + (L_{\G}+1) \sqrt{d}$.
	Hence
	\begin{align*}
		\E d_W^2 \left( \nu, \mu_{\widehat{\theta}} \right) & \leq R^2 \E d_{TV}^2 \left( \nu, \mu_{\widehat{\theta}} \right).
	\end{align*}
	The reason is: for any $f(x)$ that has Lipchitz constant $1$ with $f(0) = 0$, it is true that
	$f(x)$ is bounded in a bounded domain with radius $R$. Such a centering of $f$ is without loss of generality since $\sup_{f: f \in Lip(1)} \int f(\mu - \nu) dx = \sup_{f-f(0): f \in Lip(1)} \int (f-f(0))(\mu - \nu) dx $, for probability distributions $\mu, \nu$.
	
\end{proof}

\begin{proof}[Proof of Corollary~\ref{cor:multivariate-gaussian}]
	Suppose $\log p_\nu(x) = -\frac{1}{2}(x - b_*)' \Sigma_*^{-1} (x- b_*) + \frac{1}{2}\log \det (\Sigma_*^{-1}) - \frac{d}{2} \log (2\pi).$ And the generator class is depth-one NN, with weights $\theta = (W, b)$, $X = W Z + b$, then
	$\log p_{\mu_\theta}(x) = -\frac{1}{2} (x - b)' (W W')^{-1} (x - b) + \frac{1}{2} \log \det ((W W')^{-1}) - \frac{d}{2} \log (2\pi)$.
	
	For the discriminator, with the activation function $\sigma(t) = t^2$, one can use $O(d)$ units in a discriminator network with depth $2$, so that the two approximation error terms are zero. Note that one can also realize the quadratic activation with the ReLU activation in a bounded domain, using the construction in \cite{yarotsky2016error}. By Lemma~\ref{lem:VC-dim-bound} with degree at most $2$, $\text{VCdim}(\F) \precsim d^2 \log d$, $\text{VCdim}(\F \circ \G) \precsim (pd+d^2) \log (p + d)$. 
	Therefore $$\E d^2_{TV}\left( g_{\theta}(Z), X  \right)  \leq C \left(\frac{d^2 \log d}{n \wedge m} + \frac{(pd+d^2) \log (p + d)}{m}\right)^{1/2}.$$
\end{proof}

\subsection{Supporting Lemmas}
\label{sec:appendix}

Let's define the empirical Rademacher complexity,
\begin{align}
	\cR_n \left(\F\right) := \E_{\epsilon} \sup_{f \in \F} \frac{1}{n} \sum_{i=1}^n \epsilon_i f(X_i),
\end{align}
where $\epsilon = (\epsilon_1,\ldots,\epsilon_n)$ are i.i.d. Rademacher random variables.

\begin{lemma}[Symmetrization and entropy integral]
	\label{lem:symmetrization}
		For $\widehat{\nu}^n = \frac{1}{n}\sum_{i=1}^n \delta_{X_i}$,
		\begin{align}
			\E d_{\F}(\nu, \widehat{\nu}^n) \leq 2 \E \cR_n \left(\F\right).
		\end{align}
		Assume $\sup_{f \in \F} \| f \|_{\infty} \leq 1$, one has the standard entropy integral bound,
		\begin{align*}
			\E d_{\F}(\nu, \widehat{\nu}^n) \leq  2 \E \inf_{0<\delta<1/2} \left( 4\delta + \frac{8 \sqrt{2}}{\sqrt{n}} \int_{\delta}^{1/2}  \sqrt{\log \mathcal{N}(\epsilon, \F, \| \cdot \|_{n}) } d \epsilon \right),
		\end{align*}
		where $\| f \|_n:= \sqrt{1/n \sum_{i=1}^n f(X_i)^2}$ is the empirical $\ell_2$-metric on data $\{X_i\}_{i=1}^n$, and $\mathcal{N}(\epsilon, \F, \| \cdot \|_{n})$ is the $\ell_2$-covering number.
\end{lemma}
Remark that since $\| f \|_n \leq \max_i |f(X_i)|$, and thus $\mathcal{N}(\epsilon, \F, \| \cdot\|_n) \leq \mathcal{N}(\epsilon, \F|_{X_1 \ldots, X_n}, \infty)$. Therefore, the upper bound in the above Lemma also holds with $\mathcal{N}(\epsilon, \F|_{X_1 \ldots, X_n}, \infty)$, the $\ell_\infty$-covering number on the data.

\begin{proof}
	We use the Dudley entropy integral, a standard result in empirical process theory. For the first inequality, apply the standard symmetrization technique, we have
	\begin{align*}
		\E d_{\F}(\nu, \widehat{\nu}^n) \leq \E_{X, X'} \sup_{f \in \F} \frac{1}{n} \sum_{i=1}^n f(X_i) - f(X_i') \leq 2 \E_{X} \E_{\epsilon} \sup_{f \in \F} \frac{1}{n} \sum_{i=1}^n \epsilon_i f(X_i).
	\end{align*}
\end{proof}

The next two results, Theorems 12.2 and 14.1 in \cite{anthony2009neural}, show that the metric entropy may be bounded in terms of the pseudo-dimension.
\begin{lemma}
	\label{lem:ps-dim}
	Assume for all $f \in \F$,  $\| f \|_{\infty} \leq M$. Denote the pseudo-dimension of $\F$ as ${\rm Pdim}(\F)$, then for $n \geq {\rm Pdim}(\F)$, we have for any $\epsilon$ and any $X_1, \ldots, X_n$,
	\[
		 \mathcal{N}(\epsilon, \F|_{X_1, \ldots, X_n}, \infty) \leq \left( \frac{2e M \cdot n}{\epsilon \cdot {\rm Pdim}(\F)} \right)^{{\rm Pdim}(\F)}.
	\]
\end{lemma}

\begin{lemma}
	If $\F$ is the class of functions generated by a neural network with a fixed architecture and fixed activation functions, then
	\[
		{\rm Pdim}(\F) \leq  {\rm VCdim}(\tilde{\F})
	\]
	where $\tilde{\F}$ has only one extra input unit and one extra computation unit compared to $\F$. 
\end{lemma}

\begin{lemma}[Rademacher complexity and Pseudo-dimension]
	\label{lem:rad-vc}
	Under the condition $\max_i |f(X_i)| \leq B$, then for any $n \geq {\rm Pdim}(\F)$,
	\begin{align*}
		\cR_n(\F) \leq C \cdot B \sqrt{ \frac{{\rm Pdim}(\F) \log n}{n} }
	\end{align*}
	for some universal constant $C>0$. 
\end{lemma}
\begin{proof}
	The proof is a direct application of the Dudley entropy integral in Lemma~\ref{lem:symmetrization} and the covering number bound by pseudo-dimension in Lemma~\ref{lem:ps-dim}. See A.2.2 in \cite*{farrell2018DeepNeural} for details. 
\end{proof}

\begin{lemma}[Theorem 6 in \cite{Bartlett-etal2017_COLT}, Vapnik-Chervonenkis dimension]
	\label{lem:VC-dim-bound}
		Consider function class computed by a feed-forward neural network architecture with $W$ parameters and $U$ computation units arranged in $L$ layer. Suppose that all non-output units have piecewise-polynomial activation functions with $p+1$ pieces and degree no more than $d$, and the output unit has the identity function as its activation function. Then the VC-dimension and pseudo-dimension is upper bounded
		\[
			{\rm VCdim}(\F), {\rm Pdim}(\F) \leq C \cdot \left( LW \log (pU) + L^2 W \log d\right),
		\]
		with some universal constants $C>0$. The same result holds for pseudo-dimension ${\rm Pdim}(\F)$.	
\end{lemma}

\begin{lemma}[\cite{van2014probability}, special case of Theorem 4.8 and Example 4.9]
	\label{lem:pinsker-bobkov}
	For any two random variables $g_{\theta}(Z), X \in \mathbb{R}^d$ with $g_{\theta}(Z) \sim \mu_\theta$ and $X \sim \nu$, Pinsker's inequality asserts that
	\begin{align*}
		2 d_{TV}^2\left( \mu_\theta, \nu  \right) \leq  d_{KL} \left(  \mu_\theta || \nu  \right).
	\end{align*}
	
	Assume in addition that $Z \sim N(0, I_d)$ to be isotropic Gaussian and for all $\theta$, $\|g_\theta(z) - g_\theta(z')\| \leq L \| z - z' \|$ is $L$-Lipschitz. Then for any $X \sim \nu$ and $g_{\theta}(Z) \sim \mu_\theta$
	\begin{align*}
		d_{W}^2\left(  \nu, \mu_{\theta}  \right) \leq  2L^2 d_{KL} \left( \nu || \mu_\theta  \right).
	\end{align*}
\end{lemma}
\begin{proof}
	Consider any $1$-Lipchitz function $f: \mathbb{R}^d \rightarrow \mathbb{R}$, then $f \circ g_\theta$ is $L$-Lipschitz, which implies $f \circ g_\theta$ is $L^2$-subGaussian due to Gaussian concentration Theorem 3.25 in \cite{van2014probability}. Therefore we know $f( g_\theta(Z) )$ is $L^2$-subGaussian for any $f$ that is $1$-Lipchitz, together with Theorem 4.8 in \cite{van2014probability}, the proof completes.
\end{proof}

	\begin{lemma}[Theorem 2.5 in \cite{tsybakov2009introduction}]
		\label{lem:fano}
		Let $d(\cdot, \cdot)$ be a metric on $\Theta$.
		Assume that $H \geq 2$ and suppose $\Theta$ contains $\theta_0, \theta_1, \ldots, \theta_H$ such that:
		\begin{enumerate}
			\item $d(\theta_j, \theta_k) \geq 2s >0$, for all $j, k \in [H]$ and $j \neq k$.
			\item $\frac{1}{H} \sum_{j=1}^H d_{KL}(P_{j} || P_0) \leq c \cdot\log H$ with $0< c < 1/8$ and $P_j = P_{\theta_j}$ for $j \in [H]$.
		\end{enumerate}
		Then for any estimator $\hat{\theta}$,
		\begin{align*}
			\sup_{\theta \in \Theta} P_\theta(d(\hat{\theta}, \theta) \geq s) \geq \frac{\sqrt{H}}{1+\sqrt{H}} \left( 1 - 2c - \sqrt{\frac{2c}{\log H}} \right) > 0.
		\end{align*}
	\end{lemma}	

\end{document}